\documentclass[11pt]{amsart}
\usepackage[a4paper, footskip=1cm, headheight = 16pt, top=3.7cm, bottom=3cm,  right=2.5cm,  left=2.5cm ]{geometry}
\usepackage[utf8]{inputenc}
\usepackage{amssymb, amsthm}
\usepackage[leqno]{amsmath}
\usepackage[foot]{amsaddr}
\pdfoutput=1 
\usepackage{enumerate}
\usepackage{comment}
\usepackage{soul}
\usepackage{subcaption}
\usepackage[colorlinks,linkcolor=blue,citecolor=red]{hyperref}
\usepackage[T1]{fontenc}
\usepackage{thm-restate}

\usepackage{amsthm}
\usepackage{graphics}
\usepackage{amsmath}
\usepackage{xcolor}
\usepackage{tikz}
\usetikzlibrary{positioning}
\usetikzlibrary{snakes}

 \newtheorem{defeng}{Definition}[section]
 \newtheorem{theorem}[defeng]{Theorem}
 \newtheorem{lemma}[defeng]{Lemma}
\newtheorem{observation}[defeng]{Observation}
\newtheorem{conjecture}[defeng]{Conjecture}
\newtheorem{question}[defeng]{Question}

\DeclareMathOperator{\anchor}{anchor}
\DeclareMathOperator{\sep}{sep}

\renewcommand{\S}{\mathcal{S}}

\newcommand{\X}{\mathcal{X}}
\newcommand{\C}{\mathcal{C}}

\DeclareMathOperator{\tw}{tw}
\DeclareMathOperator{\cent}{center}

\mathchardef\mh="2D
\newcounter{tbox}

\newcommand{\sta}[1]{\vspace*{0.3cm}\refstepcounter{tbox}\noindent{ \parbox{\textwidth}{(\thetbox) \emph{#1}}}\vspace*{0.3cm}}

\newcommand{\vsp}{\vspace*{3mm}}

\usepackage{cleveref}
\crefformat{footnote}{#2\footnotemark[#1]#3}

\def\dd{\hbox{-}}

\newcommand{\mylongtitle}[1]{%
  \ifodd\value{page}%
    \protect\parbox{0.97\linewidth}{#1}\hfill%
  \else%
    \hfill\protect\parbox{0.97\linewidth}{#1}%
  \fi%
}

\title[Induced subgraphs and tree decompositions II.]{Induced subgraphs and tree decompositions\\
II. Toward walls and their line graphs in graphs of bounded degree\footnote{This is an accepted manuscript. The published paper appeared in Journal of Combinatorial Theory, Series B, Volume 164, January 2024, Pages 371--403 and is available at \url{https://doi.org/10.1016/j.jctb.2023.10.005}.}}

\author{Tara Abrishami$^{\ast \dagger}$}
\author{Maria Chudnovsky$^{\ast \dagger}$}
\author{Cemil Dibek$^{\ast \ddagger}$}
\author{Sepehr Hajebi $^{\mathsection}$}
\author{Pawe{\l} Rz\k{a}\.{z}ewski$^{\mathparagraph}$}
\author{Sophie Spirkl$^{\mathsection \parallel}$}
\author{Kristina Vu\v{s}kovi\'c $^{\ast\ast}$}
\address{$^{\ast}$Princeton University, Princeton, NJ, USA}
\address{$^{\mathsection}$Department of Combinatorics and Optimization, University of Waterloo, Waterloo, Ontario, Canada}
\address{$^{\mathparagraph}$Warsaw University of Technology, Poland/University of Warsaw, Poland. This work is a part of project
CUTACOMBS that has received funding from the European Research Council (ERC) under the European Union's
Horizon 2020 research and innovation programme (grant agreement No. 714704).}
\address{$^{\ast \ast}$School of Computing, University of Leeds, UK. Partially supported by DMS-EPSRC Grant EP/V002813/1.}
\address{$^{\dagger}$ Supported by NSF Grant DMS-1763817 and
     NSF-EPSRC Grant DMS-2120644.}
\address{$^{\ddagger}$ Supported by NSF Grant DMS-1763817.}
\address{$^{\parallel}$ We acknowledge the support of the Natural Sciences and Engineering Research Council of Canada (NSERC), [funding reference number RGPIN-2020-03912].
Cette recherche a \'et\'e financ\'ee par le Conseil de recherches en sciences naturelles et en g\'enie du Canada (CRSNG), [num\'ero de r\'ef\'erence RGPIN-2020-03912].}

\date{\today}

\allowdisplaybreaks

\begin{document}

\maketitle
\begin{abstract}

    This paper is motivated by the following question: what are the unavoidable induced subgraphs of graphs with large treewidth? Aboulker et al. made a conjecture which answers this question in graphs of bounded maximum degree, asserting that for all $k$ and $\Delta$, every graph with maximum degree at most $\Delta$ and sufficiently large treewidth contains either a subdivision of the $(k\times k)$-wall or the line graph of a subdivision of the $(k\times k)$-wall as an induced subgraph. We prove two theorems supporting this conjecture, as follows.

\begin{itemize}

\item[1.] For $t\geq 2$, a $t$\textit{-theta} is a graph consisting of two nonadjacent vertices and three internally vertex-disjoint paths between them, each of length at least $t$. A $t$\textit{-pyramid} is a graph consisting of a vertex $v$, a triangle $B$ disjoint from $v$ and three paths starting at $v$ and vertex-disjoint otherwise, each joining $v$ to a vertex of $B$, and each of length at least $t$. We prove that for all $k,t$ and $\Delta$, every graph with maximum degree at most $\Delta$ and sufficiently large treewidth contains either a $t$-theta, or a $t$-pyramid, or the line graph of a subdivision of the $(k\times k)$-wall as an induced subgraph. This affirmatively answers a question of Pilipczuk et al. asking whether every graph of bounded maximum degree and sufficiently large treewidth contains either a theta or a triangle as an induced subgraph (where a \textit{theta} means a $t$-theta for some $t\geq 2$).

\item[2.] A \textit{subcubic subdivided caterpillar} is a tree of maximum degree at most three whose all vertices of degree three lie on a path. We prove that for every $\Delta$ and subcubic subdivided caterpillar $T$, every graph with maximum degree at most $\Delta$ and sufficiently large treewidth contains either a subdivision of $T$ or the line graph of a subdivision of $T$ as an induced subgraph.

\end{itemize}

\end{abstract}




\section{Introduction} 

\begin{figure}[t!]
\centering
\includegraphics[scale=0.6]{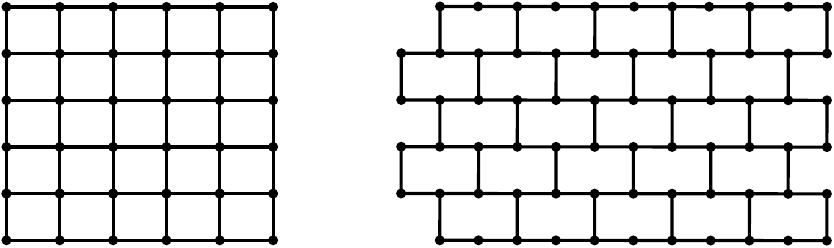}
\caption{The $6$-by-$6$ square grid (left) and the $6$-by-$6$ wall $W_{6\times 6}$ (right).}
\label{fig:Grid+Wall}
\end{figure}

All graphs in this paper are finite and simple. Let $G = (V(G), E(G))$ be a graph.  A {\em tree decomposition $(T, \beta)$} of $G$ consists of a tree $T$ and a map $\beta: V(T) \to 2^{V(G)}$, with the following properties: 
\begin{enumerate}[(i)]
    \item For every $v \in V(G)$, there exists $t \in V(T)$ such that $v \in \beta(t)$. 
    
    \item For every $v_1v_2 \in E(G)$, there exists $t \in V(T)$ such that $v_1, v_2 \in \beta(t)$.
    
    \item For every $v \in V(G)$, the subgraph of $T$ induced on the set $\beta^{-1}(v)=\{t \in V(T) \mid v \in \beta(t)\}$ is connected.
\end{enumerate}
The {\em width} of the tree decomposition $(T, \beta)$ is $\max_{v \in V(T)} |\beta(v)| -1$. The {\em treewidth} of a graph $G$, denoted by $\tw(G)$, is the minimum width of a tree decomposition of $G$. Treewidth, originally introduced by Robertson and Seymour in their study of graph minors, is widely considered to be an important graph parameter, both from a structural
\cite{RS-GMXVI} and algorithmic \cite{Bodlaender1988DynamicTreewidth} point of view. Roughly, the treewidth of a graph measures how ``close to a tree'' it is: trees have treewidth one, and in general, the larger the treewidth of a graph, the less ``tree-like'', and hence the more complicated it is. So it is natural to ask how one would certify whether a graph is of large treewidth, and in particular, what can we say about the unavoidable substructures emerging in graphs of large treewidth. As an example, for each $k$, the $k\dd by\dd k$ square grid is a planar graph of maximum degree three and with treewidth $k$, and the {\em $(k \times k)$-wall}, denoted by $W_{k \times k}$, is the $k\dd by\dd k$ hexagonal grid, which is planar graph with maximum degree three and with treewidth $k$ (the formal definition is provided at the end of Subsection~\ref{sec:defns}; see Figure \ref{fig:Grid+Wall}). Every subdivision of $W_{k \times k}$ is also a
graph of treewidth $k$. The Grid Theorem of Robertson and Seymour gives a complete characterization of the unavoidable minors, and the unavoidable subgraphs of graphs with large treewidth.
\begin{theorem}[Robertson and Seymour \cite{RS-GMV}]\label{wallminor}
For every $k\geq 1$, every graph of sufficiently large treewidth has minor isomorphic to the $k\dd by\dd k$ square grid, or equivalently, a subdivision of $W_{k \times k}$ as a subgraph.
\end{theorem}

While tree decompositions and classes of graphs with bounded treewidth are central concepts in the study of graphs with forbidden minors \cite{RS-GMXVI}, the problem of connecting tree decompositions with forbidden induced subgraphs had largely remained uninvestigated until very recently. In accordance, this work is a step toward understanding the unavoidable induced subgraphs of graphs with large treewidth. Formally, let us say a family $\mathcal{F}$ of graphs is
{\em useful} if there exists $c$ such that every graph $G$ with $\tw(G) > c$ contains a member of $\mathcal{F}$ as an induced subgraph. Then our work is motivated by the
goal of characterizing useful families. For instance, Lozin and Razgon  \cite{LR} have recently proved the following theorem, which gives a complete description of all finite useful families. Given a graph $F$, the {\em line graph} $L(F)$ of $F$ is the graph with vertex set $E(F)$, such that two vertices of $L(F)$ are adjacent if the corresponding edges of $G$ share an end.

\begin{theorem}[Lozin and Razgon \cite{LR}]\label{Lozinfinite}
  Let $\mathcal{F}$ be finite family of graphs. Then $\mathcal{F}$ is useful if and only if it contains a complete graph, a complete bipartite graph, a forest in which each component has at most three leaves, and the line graph of such a forest.
\end{theorem}

In fact, it is easy to see that the complete graph $K_t$ has treewidth $t-1$ and the complete bipartite graph $K_{t,t}$ has treewidth $t$. Also, as mentioned above, every subdivision of $W_{k \times k}$ is also of treewidth $k$, and crucially, no two non-isomorphic subdivisions of $W_{k \times k}$ are induced subgraphs of each other. The line graph of a subdivision of $W_{k \times k}$ is
another
example of a graph with large treewidth. Note that $L(W_{k \times k})$ does not contain
$W_{k \times k}$ as an induced subgraph.
In summary, if a family of graphs is useful, then it contains a complete graph, a complete bipartite graph, and for some $k$, an induced subgraph of every subdivision of $W_{k \times k}$, and an induced subgraph of the line graph of every subdivision of $W_{k \times k}$. Therefore, it would be natural to ask whether the converse of the latter statement is also true:

\begin{question}\label{usefulQ}
Let $\mathcal{F}$ be a family of graphs containing a complete graph, a complete bipartite graph, and for some $k$, an induced subgraph of every subdivision of $W_{k \times k}$, and an induced subgraph of the line graph of every subdivision of $W_{k \times k}$. Then is $\mathcal{F}$ useful?
\end{question}

It turns out that the answer to Question~\ref{usefulQ} is negative. To elaborate on this, we need a couple of definitions. By a {\em hole} in a graph we mean an induced cycle of length at least four, and an {\em even hole} is a hole on an even number of vertices. For graphs $G$ and $F$, we say that $G$ is {\em $F$-free} if $G$ does not contain an induced subgraph isomorphic to $F$. If $\mathcal{F}$ is a family of graphs, a graph $G$ is {\em $\mathcal{F}$-free} if $G$ is $F$-free for every $F \in \mathcal{F}$. 
It is not difficult to show that for large enough $k$, subdivisions of $W_{k\times k}$, line graphs of subdivisions of $W_{k\times k}$, and the complete bipartite graph $K_{k,k}$ all
contain even holes. Therefore, the following theorem provides a negative answer to Question~\ref{usefulQ}.

\begin{theorem}[Sintiari and Trotignon \cite{ST}]
\label{thm:layered_wheel}
For every integer $\ell \geq 1$, there exists an (even hole, $K_4$)-free graph $G_{\ell}$ such that $\tw(G_\ell) \geq \ell$. 
\end{theorem}

Observing that graphs $G_{\ell}$ in Theorem \ref{thm:layered_wheel} have vertices of arbitrarily large degree, the following conjecture was made (and proved for the case $\Delta\leq 3$) in \cite{Aboulker2020OnGraphs}:

\begin{conjecture}[Aboulker, Adler, Kim, Sintiari and Trotignon \cite{Aboulker2020OnGraphs}] \label{evenholestw}
  For every $\Delta > 0$ there exists $c_{\Delta}$ such that  even-hole-free graphs
  with maximum degree at most $\Delta$ have treewidth at most $c_{\Delta}$. 
\end{conjecture}

Conjecture~\ref{evenholestw} was proved in \cite{ACV} by three of
the authors of the present paper. More generally, it is conjectured in \cite{{Aboulker2020OnGraphs}} that there is an affirmative asnwer to Questoin~\ref{usefulQ} in the bounded maximum degree case (note that bounded maximum degree automatically implies that a large complete graph and a large complete bipartite graph are excluded).

\begin{conjecture}[Aboulker, Adler, Kim, Sintiari and Trotignon \cite{Aboulker2020OnGraphs}]\label{conj:wall}
For  every $\Delta>0$ there  is  a  function
$f_{\Delta}:\mathbb{N} \rightarrow \mathbb{N}$ such  that
every  graph  with  maximum degree  at  most $\Delta$ and treewidth  at  least
$f_{\Delta}(k)$
contains  a subdivision of $W_{k \times k}$ or  the line graph of a subdivision of $W_{k \times k}$ as an induced
subgraph.
\end{conjecture}

(We remark that, while the present paper was under review, Conjecture \ref{conj:wall} was proved using a different method \cite{Korhonen}. However, the techniques developed here provide foundation to a significant body of future work \cite{twx,  treealpha, twv, longclaw, c4prismfree, twiii,  twvi, twiv}.)

In \cite{Aboulker2020OnGraphs} it is proved for proper minor-closed classes of graphs (in which case the bound on the maximum degree is not needed anymore).
\begin{theorem}[Aboulker, Adler, Kim Sintiari and Trotignon \cite{Aboulker2020OnGraphs}]
For every
graph H there is a function $f_H : \mathbb{N} \to \mathbb{N}$ such that every graph of treewidth at least $f_H(k)$  and with no $H$-minor
contains a subdivision of $W_{k \times k}$  or the line graph of a subdivision of $W_{k \times k}$ as an induced subgraph.
\end{theorem}

In this paper we prove several 
theorems supporting
Conjecture~\ref{conj:wall}. In order to state our main results, we need a few more definitions.

A {\em path} is a graph $P$ with vertex set $\{p_1, \hdots, p_k\}$ and edge set $\{p_1p_2, p_2p_3, \hdots, p_{k-1}p_k\}$. We write $P = p_1 \dd \hdots \dd p_k$, and we say $p_1$ and $p_k$ are the {\em ends} of $P$. The {\em length} of the path $P$ is the number of edges in $P$. We say that $P$ is a path {\em from $p_1$ to $p_k$}, where $p_1$ and $p_k$ are the vertices of degree one in $P$. The {\em interior of $P$} is denoted $P^*$ and is defined as $P \setminus \{p_1, p_k\}$. 

Let $G$ be a graph and let $X, Y \subseteq V(G)$ be disjoint. Then, $X$ is {\em complete to $Y$} if for every $x \in X$ and $y \in Y$, we have $xy \in E(G)$, and $X$ is {\em anticomplete to $Y$} if there are no edges from $X$ to $Y$ in $G$.

The {\em claw} is the graph with vertex set $\{a, b, c, d\}$ and edge set $\{ab, ac, ad\}$. For nonnegative integers $t_1, t_2, t_3$, an $S_{t_1, t_2, t_3}$, also called a {\em long claw} or a {\em subdivided claw}, consists of a vertex $v$ and three paths $P_1, P_2, P_3$, where $P_i$ is of length $t_i$, with one end $v$, such that $V(P_1) \setminus \{v\}$, $V(P_2) \setminus \{v\}$, and $V(P_3) \setminus \{v\}$ are pairwise disjoint and anticomplete to each other. Note that for every $t$, every subdivision of $W_{k\times k}$ for large enough $k$ contains $S_{t,t,t}$ as an induced subgraph. Our first result is the following.
 
\begin{restatable}{theorem}{clawfree}
\label{thm:claw-free_nonspecific}
Let $\Delta, t, k$ be positive integers. There exists $c_{k, t, \Delta}$ such that for every $S_{t,t,t}$-free graph $G$ with maximum degree $\Delta$ and no induced subgraph isomorphic to the line graph of a subdivision of $W_{k \times k}$, we have $\tw(G) \leq c_{k, t, \Delta}$. 
\end{restatable}

A {\em theta} is a graph consisting of three internally vertex-disjoint paths $P_1 = a\dd \cdots \dd b$, $P_2 = a \dd \cdots \dd b$, and $P_3 = a \dd \cdots \dd b$ of length at least 2, such that no edges exist between the paths except the three edges incident with $a$ and the three edges incident with $b$.  A {\em $t$-theta} is a theta such that each of $P_1, P_2, P_3$ has length at least $t$. A {\em pyramid} is a graph consisting of three paths $P_1 = a \dd \cdots \dd b_1$, $P_2 = a \dd \cdots \dd b_2$, and $P_3 = a \dd \cdots \dd b_3$ of length at least 1, two of which have length at least 2, pairwise vertex-disjoint except at $a$, and such that $b_1b_2b_3$ is a triangle and no edges exist between the paths except those of the triangle and the three edges incident with $a$. A {\em $t$-pyramid} is a pyramid such that each of $P_1, P_2, P_3$ has length at least $t$. 

Note that the complete bipartite graph $K_{2,3}$ is in fact a theta. Also, for large enough $k$, every subdivision of $W_{k\times k}$ contains a theta as an induced subgraph, and the line graph of every subdivision of $W_{k\times k}$ contains a triangle. Therefore, the following theorem gives another reason why the answer to Question~\ref{usefulQ} is negative. 
\begin{theorem}[Sintiari and Trotignon \cite{ST}]
\label{thm:layered_wheel_theta}
For every integer  $\ell \geq 1$, there exists a
(theta, triangle)-free graph $G_{\ell}$ such that $\tw(G_\ell) \geq \ell$. 
\end{theorem}

In analogy to the situation with Theorem \ref{thm:layered_wheel}, the graphs $G_{\ell}$ in
Theorem~\ref{thm:layered_wheel_theta} contain vertices of arbitrary large degree. So it is asked in \cite{PSTT} whether (theta, triangle)-free graphs of bounded maximum degree have bounded treewidth (while it is proved in \cite{PSTT} that (theta, triangle,$S_{t,t,t})$-free graphs, without a bound on the maximum degree, have bounded treewidth). We give an affirmative answer to this question. Indeed, our second result, the following, establishes a far-reaching generalization of this question. It also generalizes Theorem \ref{thm:claw-free_nonspecific}, and strongly addresses Conjecture \ref{conj:wall}.

\begin{restatable}{theorem}{pyramid-theta}
\label{thm:pyramid_theta-nonspecific}
Let $\Delta, t, k$ be positive integers with $t \geq 2$. Then, there exists $c_{t, k, \Delta}$ such that for every  ($t$-theta, $t$-pyramid)-free graph $G$ with maximum degree $\Delta$ and no induced subgraph isomorphic to the line graph of a subdivision of $W_{k \times k}$, we have $\tw(G)\leq c_{t, k, \Delta}$. 
\end{restatable}

A tree $T$ is a {\em subdivided caterpillar} if there is a path $P$ in $T$ such that $P$ contains every vertex of $T$ of degree at least three in $T$. The {\em spine} of $T$ is the shortest path containing all vertices of degree at least three in $T$. A {\em leg} of a subdivided caterpillar $T$ is a path in $T$ from a vertex of degree one in $T$ to a vertex of degree at least three in $T$. A graph $G$ is {\em subcubic} if every vertex of $G$ has degree at most three. 

Note that for every subcubic subdivided caterpillar $T$ and for large enough $k$, every subdivision of $W_{k\times k}$ contains a subdivision of $T$ as an induced subgraph, and the line graph of every subdivision of $W_{k\times k}$ contains the line graph of a subdivision of $T$ as an induced subgraph. Our third result is the following.
\begin{restatable}{theorem}{caterpillar}
\label{thm:caterpillar-non_specific}
Let $\Delta$ be a positive integer and let $T$ be a subcubic subdivided caterpillar. There exists $c_{T, \Delta}$ such that for every graph $G$ with maximum degree $\Delta$ and no induced subgraph isomorphic to a subdivision of $T$ or the line graph of a subdivision of $T$, we have $\tw(G) \leq c_{T, \Delta}$. 
\end{restatable}


Let us now roughly discuss the proofs.
Usually, to prove that a certain graph family has bounded treewidth, one attempts to construct a collection of ``non-crossing decompositions,'' which  roughly  means that the decompositions
``cooperate'' with each other, and the pieces that are obtained when the graph
is simultaneously decomposed by all the decompositions in the collection
``line up'' to form a tree structure. Such collections of decompositions are
called ``laminar.''
In all the cases above, there is a natural family of decompositions to turn to,
sharing a certain structural property: all the decompositions arise from
removing from the graph the  neighborhood of a small connected subgraph.
Unfortunately,
these natural decompositions are very far from being non-crossing, and therefore
they cannot be used in traditional ways to get tree-decompositions. What turns out to be true, however, is that, due to the bound on the maximum degree of the graph, these collections of decompositions can be partitioned into a bounded number
of laminar collections (where the bound on the number of collections depends on the maximum degree and on the precise nature of the decomposition). We will explain how to make use of this fact in Section  \ref{sec:central_bags}.
\subsection*{Structure of the paper}
We begin in Section \ref{sec:defns} with a review of relevant definitions and notation. In Section \ref{sec:balanced-separators}, we define an important graph parameter tied to treewidth called separation number. In Section \ref{sec:central_bags} we prove Theorem~\ref{thm:centralbag}, which summarizes our main proof method. In Section \ref{sec:tw_of_clawfree}, we bound the treewidth of claw-free graphs with  no line graph of a subdivision of a wall, and in Section \ref{sec:claw_free_result}, we apply the results of Section \ref{sec:tw_of_clawfree} to prove Theorem \ref{thm:claw-free_nonspecific}. In Section \ref{sec:theta_pyramid}, we prove Theorem \ref{thm:pyramid_theta-nonspecific}, and in Section \ref{sec:caterpillar}, we prove Theorem \ref{thm:caterpillar-non_specific}. 

\subsection{Definitions and Notation}
\label{sec:defns}
Let $G$ be a graph. In this paper, we use vertex sets and their induced subgraphs interchangeably. Let $H$ be a graph. We say that $X \subseteq V(G)$ {\em is an $H$ in $G$} if $X$ is isomorphic to $H$. We say that $G$ {\em contains} $H$ if there exists $X \subseteq V(G)$ such that $X$ is an $H$ in $G$. 

The {\em open neighborhood} of a vertex $v \in V(G)$, denoted $N(v)$, is the set of all vertices adjacent to $v$. The {\em degree} of $v \in V(G)$ is the size of its open neighborhood. A graph $G$ has {\em maximum degree $\Delta$} if the degree of every vertex $v \in V(G)$ is at most $\Delta$. The {\em closed neighborhood} of a vertex $v \in V(G)$ is denoted $N[v]$ and is defined as $N[v] = N(v) \cup \{v\}$. Let $X \subseteq V(G)$. The {\em open neighborhood of $X$}, denoted $N(X)$, is the set of all vertices of $G \setminus X$ with a neighbor in $X$. The {\em closed neighborhood of $X$} is denoted $N[X]$ and is defined as $N[X] = N(X) \cup X$.  

A set $X \subseteq V(G)$ is {\em connected} if for every $x, y \in X$, there is a path $P$ in $X$ from $x$ to $y$. A set $C \subseteq V(G)$ is a {\em cutset} of a connected graph $G$ if $G \setminus C$ is not connected. A set $D$ is a {\em connected component of $G$} if $D$ is inclusion-wise maximal such that $D \subseteq V(G)$ and $D$ is connected. 

Let $u, v \in V(G)$ and let $X \subseteq V(G)$. The {\em distance between $u$ and $v$} is the length of a shortest path from $u$ to $v$ in $G$. The {\em distance between $u$ and $X$} is the length of a shortest path from $u$ to a vertex $x \in X$ in $G$. We denote by $N^d(v)$ the set of vertices at distance exactly $d$ from $v$ in $G$, and by $N^d[v]$ the set of vertices at distance at most $d$ from $v$ in $G$. Similarly, we denote by $N^d[X]$ the set of vertices of distance at most $d$ from $X$ in $G$. The {\em diameter} of a connected set $X \subseteq V(G)$ is the maximum distance in $G$ between two vertices of $X$.


A {\em clique} is a set $K \subseteq V(G)$ such that every pair of vertices in $K$ is adjacent. An {\em independent set} is a set $I \subseteq V(G)$ such that every pair of vertices in $I$ is non-adjacent. The {\em clique number of $G$}, denoted $\omega(G)$, is the size of a largest clique in $G$. The {\em independence number of $G$}, denoted $\alpha(G)$, is the size of a largest independent set in $G$. 

A {\em weight function on $G$} is a function $w:V(G) \to \mathbb{R}$ that assigns a non-negative real number to every vertex of $G$. A weight function is {\em normal} if $w(V(G)) = 1$. Unless otherwise specified, we assume all weight functions are normal. We denote by $w^{\max}$ the maximum weight of a vertex; i.e. $w^{\max} = \max_{v \in V(G)} w(v)$. 

Finally, let us include the precise definition of a wall. The {\em $(n \times m)$-wall}, denoted $W_{n \times m}$, is the graph $G$ with vertex set
\begin{align*}
    V(G) =&
\{(1, 2j - 1) \mid 1 \leq j \leq m \} \\ 
&\cup \{(i, j) \mid 1 < i < n, 1 \leq j \leq 2m\} \\
&\cup \{(n, 2j - 1) \mid 1 \leq j \leq m, \text{ if $n$ is even}\} \\	& \cup \{(n, 2j) \mid 1 \leq j \leq m, \text{ if $n$ is odd }\}
\end{align*}

and edge set
\begin{align*}
    E(G) =& \{(1, 2j - 1),(1, 2j + 1) \mid 1 \leq j \leq m - 1\}\\
& \cup \{(i, j),(i, j + 1) \mid 2 \leq i < n, 1 \leq j < 2m \} \\
&\cup \{(n, 2j),(n, 2j + 2)) \mid 1\leq j < m \text{ if $n$ is odd} \} \\	
&\cup \{(n, 2j - 1),(n, 2j + 1) \mid 1 \leq j < m \text{ if $n$ is even} \} \\	
&\cup \{(i, j),(i + 1, j) \mid 1 \leq i < n, 1 \leq j \leq 2m, i, j \text{ odd} \} \\
&\cup \{(i, j),(i + 1, j) \mid 1 \leq i < n, 1 \leq j \leq 2m, i, j \text{ even} \}.
\end{align*}
Again, see Figure \ref{fig:Grid+Wall} for an example.
\subsection{Balanced separators and treewidth}
\label{sec:balanced-separators}
Treewidth is tied to a parameter called the separation number. Let $G$ be a graph, let $S \subseteq V(G)$, let $k$ be a positive integer, and let $c \in [\frac{1}{2}, 1)$. A set $X \subseteq V(G)$ is a {\em $(k, S, c)^*$-separator} if $|X| \leq k$ and for every component $D$ of $G \setminus X$, it holds that $|D \cap S| \leq c|S|$. The {\em separation number} $\sep_c^*(G)$ is the minimum $k$ such that $G$ has a $(k, S, c)^*$-separator for every $S \subseteq V(G)$. The following lemma states that the separation number gives an upper bound for the treewidth of a graph. 

\begin{lemma}[Harvey and Wood \cite{HarveyWood}]
\label{lemma:harvey-wood}
For every $c \in [\frac{1}{2}, 1)$ and every graph $G$, we have $\tw(G) + 1 \leq \frac{1}{1-c} \sep_c^*(G)$. 
\end{lemma}

Now, we redefine $(k, S, c)^*$-separators using weight functions.  Given a normal weight function $w$ on a graph $G$ and a constant $c \in [\frac{1}{2}, 1)$, a set $X \subseteq V(G)$ is a {\em $(w, c)$-balanced separator of $G$} if $w(D) \leq c$ for every component $D$ of $G \setminus X$.

  We call a weight function $w$ on $G$ a {\em uniform weight function} if there exists $Y \subseteq V(G)$ such that $w(v) = \frac{1}{|Y|}$ if $v \in Y$, and $w(v) = 0$ if $v \not \in Y$. Lemma \ref{lemma:harvey-wood} implies the following: 

\begin{lemma}
\label{lemma:harvey-wood-weights}
Let $c \in [\frac{1}{2}, 1)$ and let $G$ be a graph. If $G$ has a $(w, c)$-balanced separator of size at most $k$ for every uniform weight function $w$, then $\tw(G) \leq \frac{1}{1-c}k$. 
\end{lemma}
\begin{proof}
 We prove that $\sep_c^*(G) \leq k$. Let $S \subseteq V(G)$ and let $w_S$ be the weight function on $G$ such that $w_S(v) = \frac{1}{|S|}$ if $v \in S$, and $w_S(v) = 0$ otherwise. Since $w_S$ is a uniform weight function, it follows that $G$ has a $(w_S, c)$-balanced separator $X$ such that $|X| \leq k$. Let $D$ be a component of $G \setminus X$, so $w(D) \leq c$. Consequently, $|D \cap S| \leq c|S|$, and so $X$ is a $(k, S, c)^*$-separator. Therefore, $\sep_c^*(G) \leq k$, and the result follows from Lemma \ref{lemma:harvey-wood}.
\end{proof}
Lemma \ref{lemma:harvey-wood-weights} implies that if for some fixed $c \in [\frac{1}{2}, 1)$, $G$ has a balanced separator of size $k$ for every weight function $w$, then the treewidth of $G$ is bounded by a function of $k$. The next lemma states the converse. 

\begin{lemma}[Cygan, Fomin, Kowalik, Lokshtanov, Marx, Pilipczuk and Pilipczuk \cite{PA}]
\label{lemma:bounded-tw-balanced-separator}
    If $\tw(G) \leq k$, then $G$ has a $(w, c)$-balanced separator of size at most $k+1$ for every normal weight function $w$ and for every $c \in [\frac{1}{2}, 1)$.
\end{lemma}
Together, Lemmas \ref{lemma:harvey-wood-weights} and \ref{lemma:bounded-tw-balanced-separator} show that treewidth is tied to the size of balanced separators. In this paper, we rely on balanced separators to prove that graphs have bounded treewidth. In what follows, we will often assume that $G$ has no $(w, c)$-balanced separator of size $d$ for some normal weight function $w$, $c \in [\frac{1}{2}, 1)$, and positive integer $d$, since otherwise, in light of Lemma \ref{lemma:harvey-wood-weights}, we are done. 



\section{Central bags and forcers}
\label{sec:central_bags}
A {\em separation} of a graph $G$ is a triple $(A, C, B)$ with $A, C, B \subseteq V(G)$ such that $A$, $C$, and $B$ are pairwise disjoint, $A \cup C \cup B = V(G)$, and $A$ is anticomplete to $B$. If $A$ and $B$ are both non-empty, then $C$ is a cutset of $G$. If $S = (A, C, B)$ is a separation, we write $A(S) = A$, $C(S) = C$, and $B(S) = B$.

Separations provide a way to organize the structure of connected graphs around important cutsets. To that end, it is often useful to characterize the relationship between two separations of a graph. Here, we define two relations between graph separations. Two separations $S_1 = (A_1, C_1, B_1)$ and $S_2 = (A_2, C_2, B_2)$ of $G$ are {\em non-crossing} if (possibly exchanging the roles of
$A_1$ and $B_1$, and of $A_2$ and $B_2$) $A_1 \cap C_2 = \emptyset$, $A_2 \cap C_1 = \emptyset$, and $A_1 \cap A_2 = \emptyset$; and  {\em loosely non-crossing} if $A_1 \cap C_2 = \emptyset$ and $A_2 \cap C_1 = \emptyset$.
The two non-crossing properties are illustrated in Figure \ref{fig:non-crossing}. Note that if two separations are non-crossing, then they are also loosely non-crossing.
A collection $\S$ of separations is {\em (loosely) laminar} if the separations of $\S$ are pairwise (loosely) non-crossing. 

\begin{figure}[ht]
\begin{subfigure}{0.48\textwidth}   
\centering 
\begin{tabular}{c|c|c|c}
        & $A_1$ & $C_1$ & $B_1$ \\ \hline
    $A_2$    & $\varnothing$ & $\varnothing$ & \\ \hline
    $C_2$ & $\varnothing$ & & \\ \hline
    $B_2$ & & & \\
    \end{tabular}
    \caption{Non-crossing}
\end{subfigure}\quad 
\begin{subfigure}{0.48\textwidth}
\centering
        \begin{tabular}{c|c|c|c}
        & $A_1$ & $C_1$ & $B_1$ \\ \hline
    $A_2$    & & $\varnothing$ & \\ \hline
    $C_2$ & $\varnothing$ & & \\ \hline
    $B_2$ & & &\\
    \end{tabular}
    \caption{Loosely non-crossing}
\end{subfigure}
\caption{Illustrations of two separations $S_1 = (A_1, C_1, B_1)$ and $S_2 = (A_2, C_2, B_2)$ being (A) non-crossing and (B) loosely non-crossing.
}
    \label{fig:non-crossing}
\end{figure}

The notion of non-crossing separations was used  in \cite{RS-GMX} as part of the study of tree decompositions. Let $(T, \beta)$ be a tree decomposition of $G$. For $X \subseteq V(T)$, let $\beta(X) = \bigcup_{x \in X} \beta(x)$. Let $t_1t_2 \in E(T)$ be an edge of $T$. Let $T_1$ and $T_2$ be the components of $T \setminus \{t_1t_2\}$ containing $t_1$ and $t_2$, respectively. Let $C = \beta(t_1) \cap \beta(t_2)$, let $A = \beta(V(T_1)) \setminus C$, and let $B = \beta(V(T_2)) \setminus C$. Then, it follows from the properties of tree decompositions that $(A, C, B)$ is a separation of $G$. Up to symmetry between $A$ and $B$, we say that $(A, C, B)$ is the {\em separation of $G$ corresponding to the edge $t_1t_2$ of $T$}. Let $\S(T, \beta)$ be the collection of separations corresponding to the edges of $T$. It is shown in  \cite{RS-GMX} that for every tree decomposition $(T, \beta)$ of $G$, it holds that $\S(T, \beta)$ is laminar,  and conversely, for every laminar collection $\S$ of separations of $G$, there exists a tree decomposition $(T, \beta)$ of $G$ such that $\S = \S(T, \beta)$. Therefore, there is a one-to-one correspondence between the laminar collections of separations of a graph $G$ and the tree decompositions of $G$. This correspondence can be used to show that a graph has bounded treewidth: instead of attempting to bound the treewidth directly, one can instead study
its separations.

In this paper we modify the approach above and use loosely laminar
collections of separations (in fact, a slight variant of that). We then introduce a tool that reduces the problem of
bounding the treewidth of a graph $G$ to the problem of bounding the treewidth of a certain induced subgraph $\beta$ of $G$. This is summarized in
Theorem~\ref{thm:centralbag}, but we try to explain and motivate our lemmas and definitions as we go.

Let $G$ be a graph and let $w:V(G) \rightarrow \mathbb{R}$ be a normal weight function on $G$.
 Let $\varepsilon \in (0, \frac{1}{2}]$. We say that a separation $(A, C, B)$ is
 {\em $\varepsilon$-skewed} if $w(A) < \varepsilon$ or if $w(B) < \varepsilon$.
 For the remainder of the paper, we assume by convention that if $S$ is
 $\varepsilon$-skewed for some $\varepsilon \in (0, \frac{1}{1}]$, then $w(A(S)) < \varepsilon$. We have now broken
 the symmetry between $A$ and $B$. Let us say that
 separations $(A_1, C_1, B_1)$ and $(A_2,C_2,B_2)$ are {\em $A$-loosely non-crossing}
 if $A_1 \cap C_2=A_2 \cap C_1=\emptyset$ and {\em $A$-non-crossing} if
  $A_1 \cap C_2=A_2 \cap C_1=A_1 \cap A_2=\emptyset$ (we also sometimes use ``$A$-loosely crossing'' and ``$A$-crossing'' to mean  ``not $A$-loosely non-crossing'' and ``not $A$-non-crossing.'')
  We define $A$-laminar and $A$-loosely laminar collections of separations
  similarly.
 Note that if $(A_1, C_1, B_1)$ and $(A_2, C_2, B_2)$ are $A$-loosely non-crossing $\varepsilon$-skewed separations, then every component $D$ of $A_1 \cup A_2$ is a component of $A_1$ or a component of $A_2$, and therefore $w(D) < \varepsilon$.

Let $G$ be a graph, $w$ a normal weight function on $G$,  $c \in [\frac{1}{2}, 1)$ and $d$ an integer such that $G$ has no $(w,c)$-balanced separator of size at most $d$ (this will be our standard set up, in view of
Lemma~\ref{lemma:harvey-wood-weights}).
Let $S=(A,C,B)$ be a separation of $G$ with $|C| \leq d$. Since $C$ is not a $(w, c)$-balanced separator of $G$, we deduce that  $w(A) > c$ or $w(B) > c$.  Therefore, $S$ is $(1-c)$-skewed, so by our convention $w(A) < 1-c$, and consequently $w(B)>c$.

Henceforth, we assume that all collections of separations are ordered, and we refer to these ordered collections as {\em sequences of separations}. Let $G$ be a connected graph, let $w$ be a normal weight function on $G$, and let $\S$ be
an  $A$-loosely laminar sequence of $\varepsilon$-skewed separations of $G$. The {\em central bag for $\S$}, denoted $\beta_\S$, is defined as follows: 
$$\beta_\S = \bigcap_{S \in \S} B(S) \cup C(S).$$


We also define a weight function $w_\S:\beta_\S \to [0, 1]$ on $\beta_\S$ as follows. Let $\S = (S_1, \hdots, S_k)$, and assume that $C(S_i) \neq \emptyset$ for all $S_i \in \S$. We equip the sequence $\S$ with an {\em anchor map} $\anchor_\S: \S \to V(G)$, such that $\anchor_\S(S_i) \in C(S_i)$ for every $S_i \in \S$. We call $\anchor_\S(S)$ the {\em anchor for $S$}. Since $\S$ is loosely laminar, it holds that $\anchor_\S(S) \in \beta_\S$ for every $S \in \S$ (this is proven in \eqref{C-connected} of Lemma \ref{lemma:central_bag_1}).  For $v \in \beta_\S$, let $a(v) = \{i \ \text{s.t.} \ v \text{ is the anchor for } S_i\}$. Let $w^*(A_i) = w(A_i \setminus \bigcup_{1 \leq j < i}A_j)$. Then, we let $w_\S(v) = w(v) + \sum_{i \in a(v)} w^*(A_i)$ for all $v \in \beta_\S$. Thus, the anchor for a separation $S$ is a way to record the weight of $A(S)$ in $\beta_\S$. 

We state the next few lemmas in slightly greater generality than what we need here, in order to be able to use them in future work.
The following lemma gives important properties of central bags and their weight functions. 

\begin{lemma} 
Let $c \in [\frac{1}{2}, 1)$ and let $ d$ be a positive integer. Let $G$ be a connected graph, let $w$ be a normal weight function on $G$, and suppose $G$ has no $(w, c)$-balanced separator of size at most $d$. Let $\mathcal{S}$ be an $A$-loosely laminar sequence of separations of $G$ such that $C(S)$ is connected and $|C(S)| \leq d$ for all $S \in \S$, and let $\beta_\S$ be the central bag for $\S$. Then,
\begin{enumerate}[(i)]
\itemsep -0.2em
\item \label{C-connected} $C(S) \subseteq \beta_\S$ for every $S \in \S$,

\item \label{betaS-connected} $\beta_\S$ is connected,

\item \label{normal-weight} $w_\S(\beta_\S) = 1$.
\end{enumerate}
\label{lemma:central_bag_1}
\end{lemma}
\begin{proof}
Let $\S = (S_1, \hdots, S_k)$, and let $S_i = (A_i, C_i, B_i)$ for all $1 \leq i \leq k$. Since $\S$ is a $A$-loosely laminar sequence of separations, we have $C_i \cap A_j = \emptyset$ for all $1 \leq i, j \leq k$. Since $V(G) \setminus \beta_\S \subseteq \bigcup_{1 \leq i \leq k} A_i$, it follows that $C_i \cap (V(G) \setminus \beta_\S) = \emptyset$. Therefore, $C_i \subseteq \beta_\S$ for all $1 \leq i \leq k$. This proves \eqref{C-connected}. 

Let $D$ be a connected component of $\beta_\S$. Let $I = \{i: C_i \cap D \neq \emptyset\}$. Since $C_i \subseteq \beta_\S$ and $C_i$ is connected, it follows that $C_i \subseteq D$ for all $i \in I$. Since $N(A_i) \subseteq C_i$, we deduce that $D \cup \bigcup_{i \in I} A_i$ contains a connected component of $G$. Since $G$ is connected, we have $D \cup \bigcup_{i \in I} A_i = V(G)$, and so $D = \beta_\S$. This proves \eqref{betaS-connected}. 

For $1 \leq i \leq k$, let $v_i = \anchor_\S(S_i)$. Note that
\begin{align*}
w_\S(\beta_\S)
&= \sum_{v \in \beta_\S \setminus \{v_1, \hdots, v_k\}} w_\S(v) + \sum_{v \in \{v_1, \hdots, v_k\}} w_\S(v) \\
&= \sum_{v \in \beta_\S} w(v) + \sum_{1 \leq i \leq k} w\left(A_i \setminus \bigcup_{1 \leq j < i} A_j\right) \\
&= \sum_{v \in V(G)} w(v),
\end{align*}
where the last equality holds since for all $v \not \in \beta_\S$, we have $v \in A_i$ for some $1 \leq i \leq k$. Since $w(G) = 1$, it follows that $w_\S(\beta_\S) = 1$. This proves \eqref{normal-weight}.
\end{proof}

\begin{lemma}
\label{lemma:central_bag_2}
Let $c \in [\frac{1}{2}, 1)$, and let $\Delta, d, t$ be positive integers with $d \geq 1 + \Delta + \hdots + \Delta^{t}$. Let $G$ be a connected graph with maximum degree $\Delta$, let $w$ be a normal weight function on $G$, and suppose $G$ has no $(w, c)$-balanced separator of size at most $d$. Let $\S$ be an $A$-loosely laminar sequence of separations such that for every $S \in \S$, it holds that $C(S)$ is connected and has diameter at most $t$.
Let $\beta_S$ be the central bag for $\S$. Then, $\beta_\S$ has no $(w_\S, c)$-balanced separator of size at most $d(1 + \Delta + \hdots + \Delta^{t})^{-1}$.
\end{lemma}
\begin{proof}
Suppose for a contradiction that $\beta_\S$ has a $(w_\S, c)$-balanced separator $Y$ of size at most $d(1 + \Delta + \hdots + \Delta^t)^{-1}$. Since $G$ has no $(w, c)$-balanced separator of size at most $d$, it follows that $N^t[Y]$ is not a $(w, c)$-balanced separator of $G$ of size at most $d$. Since $|N^t[Y]| \leq |Y|(1 + \Delta + \hdots + \Delta^t) \leq d$, it follows that there exists a connected component $X$ of $G \setminus N^t[Y]$ such that $w(X) > c$. Let $Q_1, \hdots, Q_\ell$ be the connected components of $\beta_\S \setminus Y$. Let $D_1, \hdots, D_m$ be the connected components of $G \setminus \beta_\S$. Let $\mathcal{I} = \{i \text{ s.t. } Q_i \cap X \neq \emptyset\}$ and $\mathcal{J} = \{j \text{ s.t. } D_j \cap X \neq \emptyset\}$. Since $\S$ is $A$-loosely laminar, it holds that for every $1 \leq j \leq m$, there exists $S \in \S$ such that $D_j \subseteq A(S)$. Let $\S = (S_1, \hdots, S_k)$. For $j \in \mathcal{J}$, let $f(j)$ be minimum such that $D_j \subseteq A(S_{f(j)})$, let $S(j) = S_{f(j)}$, and let $v(j) = \anchor_\S(S(j))$. 

\sta{\label{C_jcapXempty} For all $j \in \mathcal{J}$, it holds that $C(S(j)) \cap Y = \emptyset$.}

Suppose $C(S(j)) \cap Y \neq \emptyset$ for some $j \in \mathcal{J}$. By assumption, $C(S(j))$ has diameter at most $t$, so $C(S(j)) \subseteq N^t[Y]$. Then, $N(D_j) \subseteq N^t[Y]$, and so $D_j = X$. By the choice of $d$ and since $D_j \subseteq A(S(j))$, it follows that $w(D_j) \leq w(A(S(j))) < 1-c \leq c$, a contradiction. This proves \eqref{C_jcapXempty}.  \\

Suppose that $|\mathcal{I}| = 0$. It follows that $|\mathcal{J}| \neq 0$, so let $j \in J$. Now, $D_j \subseteq A(S(j))$ and, since $\mathcal{I} = \emptyset$, it holds that  $C(S(j)) \subseteq Y$, contradicting \eqref{C_jcapXempty}. Now, suppose that $|\mathcal{I}| \geq 2$. Assume $Q_1, Q_2$ are such that $Q_1 \cap X \neq \emptyset$ and $Q_2 \cap X \neq \emptyset$. Since $Q_1$ and $Q_2$ are distinct connected components of $\beta_\S \setminus Y$, it follows that there exists $j \in \mathcal{J}$ such that $N(D_j) \cap Q_1 \neq \emptyset$ and $N(D_j) \cap Q_2 \neq \emptyset$. Since $N(D_j) \subseteq C(S(j))$, it follows that $C(S(j)) \cap Q_1 \neq \emptyset$ and $C(S(j)) \cap Q_2 \neq \emptyset$. Since $C(S(j))$ is connected, it holds that $C(S(j)) \cap Y \neq \emptyset$, contradicting \eqref{C_jcapXempty}. 

Therefore, $|\mathcal{I}| = 1$. Let $\mathcal{I} = \{i\}$. It follows that $C(S(j)) \subseteq Q_i$ for all $j \in \mathcal{J}$, and, in particular, $v(j) \in Q_i$ for all $j \in \mathcal{J}$. Let $a(v) = \{t \text{ s.t. } v = \anchor_\S(S_t)\}$. Now, 

\begin{align*} w_\S(Q_i) &= w(Q_i) + \sum_{v \in Q_i} \sum_{t \in a(v)} w^*(A(S_t)) \\
& \geq w(Q_i) + \sum_{j \in \mathcal{J}} w^*(A(S(j))) \\
& \geq w(Q_i) + \sum_{j \in \mathcal{J}} w(D_j) \geq w(X).
\end{align*}
But $w_\S(Q_i) \leq c$, since $Y$ is a $(w_\S, c)$-balanced separator of $\beta_\S$, and so $w(X) \leq c$, a contradiction. This proves Lemma \ref{lemma:central_bag_2}.
\end{proof}

Lemma \ref{lemma:central_bag_2} shows that if $G$ is a graph with large treewidth, then the central bag $\beta_\S$ for a well-behaved sequence $\S$ of separations of $G$ also has large treewidth. Our goal is to construct a sequence of separations $\S$ such that bounding the treewidth of $\beta_\S$ is easier than bounding the treewidth of $G$. We discuss how to find such a sequence later in this section. However, our candidate sequences of separations are usually not $A$-loosely laminar. Therefore, we would like a generalization of Lemma \ref{lemma:central_bag_2} that holds for sequences of separations that are not necessarily $A$-loosely laminar. We do this by defining the dimension of sequences of separations. Let $G$ be a graph and let $\S$ be a sequence of separations of $G$. The {\em dimension of $\S$}, denoted $\dim(\S)$, is the minimum number of laminar sequences of separations with union $\S$. Clearly, $\dim(\S) = 1$ if and only if $\S$ is laminar.

Next, we need the notion of a canonical separation. Let $G$ be a graph, let $X \subseteq V(G)$, and let us fix an ordering $(v_1, \hdots, v_n)$ of $V(G)$. The ordering defines a lexicographic order on the subsets of $V(G)$.
The {\em canonical separation for $X$}, denoted $S_X = (A_X, C_X, B_X)$, is defined as follows: $B_X$ is the largest-weight connected component of $G \setminus N[X]$, $C_X = X \cup (N[X] \cap N(B))$, and $A_X = V(G) \setminus (B_X \cup C_X)$. If there is more than one largest weight connected component of $G \setminus N[X]$, we choose the lexicographically minimum largest-weight component. Note that the definition of canonical separation is compatible with the convention that if a separation $S$ is $\varepsilon$-skewed, then $w(A) < \varepsilon$. The set $X$ is called the {\em center} of the separation $S_X$ and is denoted $\cent(S_X)$. A separation $S$ is called a {\em canonical separation} if there exists $X \subseteq V(G)$ such that $S = S_X$. For the remainder of the paper, if $S_X$ is a canonical separation, then we assume that the anchor for $S_X$ is contained in $\cent(S_X)$. 

Let $\S$ be an $A$-loosely laminar sequence of canonical separations, and let $S_1, S_2 \in \S$. Suppose $B(S_1) \cup C(S_1) \subseteq B(S_2) \cup C(S_2)$. Then, $\beta_{\S \setminus S_2} = \beta_\S$, so $S_2$ is an ``unnecessary'' member of $\S$ with respect to the central bag.  We say that $S_1$ is a {\em shield for $S_2$} if $B(S_1) \cup C(S_1) \subseteq B(S_2) \cup C(S_2)$.  

\begin{lemma}
Let $c \in [\frac{1}{2}, 1)$ and let $\Delta, d, t$ be positive integers with $d \geq 1 + \Delta + \hdots + \Delta^{t+1}$. Let $G$ be a connected graph with maximum degree $\Delta$, let $w$ be a normal weight function on $G$, and suppose $G$ has no $(w, c)$-balanced separator of size at most $d$. Let $S_X = (A_X, C_X, B_X)$ be a canonical separation of $G$ such that $C_X$ has diameter at most $t$. Then, $w(B_X) > c$. 
\label{lemma:B_large}
\end{lemma}
\begin{proof}
Since $X \subseteq C_X \subseteq N[X]$ and $C_X$ has diameter at most $t$, it follows that $N[X]$ has diameter at most $t+1$, and so $|N[X]| \leq 1 + \Delta + \hdots + \Delta^{t+1} \leq d$. Suppose $w(B_X) \leq c$. By definition, $B_X$ is a largest-weight connected component of $G \setminus N[X]$, so $w(D) \leq c$ for every connected component $D$ of $G \setminus N[X]$. But now $N[X]$ is a $(w, c)$-balanced separator of $G$ of size at most $d$, a contradiction. This proves the lemma. 
\end{proof}

Let $\S$ be a sequence of separations. We say $\S$ is {\em primordial} if for every distinct $S_1, S_2 \in \S$, it holds that $S_1$ is not a shield for $S_2$. Note that $S_1$ is a shield for $S_1$. Also, we say  $\S$ is {\em $(a, t)$-good} if every vertex $v \in V(G)$ is the anchor for at most $a$ separations in $\S$ and if $C(S)$ has diameter at most $t$ for every $S \in \S$.

\begin{lemma}
Let $c \in [\frac{1}{2}, 1)$ and let $a, \Delta, d, t$ be positive integers with $d \geq 1 + \Delta + \hdots + \Delta^{t+1}$. Let $G$ be a connected graph with maximum degree $\Delta$, let $w$ be a normal weight function on $G$, and suppose $G$ has no $(w, c)$-balanced separator of size at most $d$. Let $\S$ be an $(a, t)$-good primordial laminar sequence of canonical separations. Then, $\S$ is $A$-laminar. 
\label{lemma:primordial-to-Alaminar}
\end{lemma}
\begin{proof}
Suppose $S, S' \in \S$ are such that $S, S'$ are $A$-crossing. By Lemma \ref{lemma:B_large}, it holds that $w(B(S)) > c$ and $w(B(S')) > c$, so it follows that $B(S) \cap B(S') \neq \emptyset$. Therefore, since $S, S'$ are non-crossing but $A$-crossing, we may assume up to symmetry between $S$ and $S'$ that $A(S') \cap C(S) = A(S') \cap B(S) = B(S) \cap C(S') = \emptyset$. But now $B(S) \cup C(S) \subseteq B(S') \cup C(S')$, so $S$ is a shield for $S'$, a contradiction.  
\end{proof}


The following lemma extends the idea of a central bag to sequences of separations of bounded dimension. 

\begin{lemma}
\label{lemma:central_bag_high_dimension}
Let $c \in [\frac{1}{2}, 1)$ and let $d, \Delta, k, a, t$ be positive integers with $d \geq (1 + \Delta + \hdots + \Delta^{t+1})(1 + \Delta + \hdots + \Delta^{t})^{k}$. Let $G$ be a connected graph with maximum degree $\Delta$, let $w$ be a normal weight function on $G$, and suppose $G$ has no $(w, c)$-balanced separator of size at most $d$. Let $\S$ be an $(a, t)$-good sequence of canonical separations of $G$ with $\dim(\S) = k$. Assume that $\cent(S)$ is connected for every $S \in \S$. Then, there exists a sequence $\S_1, \hdots, \S_k$ of $A$-laminar sequences of separations, with $\S^* = \S_1 \cup \hdots \cup \S_k$ and $\beta = \bigcap_{S \in \S^*} B(S) \cup C(S)$, such that the following hold: 
\begin{enumerate}[(i)]
    \item \label{S-star-in-S} $\S^* \subseteq \S$,
    
    
    
    \item \label{S-star-is-covering} For all $S \in \S \setminus \S^*$, there exists $S' \in \S^*$ such that either $S'$ is a shield for $S$ or $\cent(S) \cap A(S') \neq \emptyset$,
    
    \item \label{beta-connected} $\beta$ is connected,
    
    \item \label{beta-no-balanced-separator} There is a normal weight function $w_\S$ on $\beta$ such that $\beta$ has no $(w_\S, c)$-balanced separator of size at most $d(1 + \Delta + \hdots + \Delta^{t})^{-k}$. 
\end{enumerate}
\end{lemma}
\begin{proof}
Let $\S_1', \hdots, \S_k'$ be a partition of $\S$ into laminar sequences. First, we will prove inductively that there exists an $A$-laminar sequence $\S_i$ for every $i \in \{1, \hdots, k\}$, with $\S_i \subseteq \S_i'$, and $\beta_i = \bigcap_{S \in \S_1 \cup \hdots \cup \S_i} B(S) \cup C(S)$,
such that $\beta_i$ is connected and there exists a normal weight function $w_i$ on $\beta_i$ such that $\beta_i$ has no $(w_i, c)$-balanced separator of size at most $d(1 + \Delta + \hdots + \Delta^{t+1})^{-i}$. This will prove \eqref{beta-connected} and \eqref{beta-no-balanced-separator}. 

Let $T_1 = \{B(S) \cup C(S) : S \in \S_1'\}$, and let $\S_1$ be the sequence formed by adding, for every inclusion-wise minimal $Y \in T_1$, a separation $S \in \S_1'$ such that $B(S) \cup C(S) = Y$. Notice that by construction, no separation in $\S_1$ is a shield for another separation in $\S_1$, so $\S_1$ is primordial. By Lemma \ref{lemma:primordial-to-Alaminar}, $\S_1$ is $A$-laminar. Let $\beta_1$ be the central bag for $\S_1$ (so $\beta_1 = \bigcap_{S \in \S_1} B(S) \cup C(S)$), and let $w_1$ be the weight function $w_{\S_1}$ on $\beta_1$. Note that since $\S$ is $(a, t)$-good and $d \geq 1 + \Delta + \hdots + \Delta^{t+1}$, it follows that $C(S) \leq d$ for all $S \in \S$. By Lemmas \ref{lemma:central_bag_1} and \ref{lemma:central_bag_2}, $\beta_1$ is connected and $\beta_1$ has no $(w_1, c)$-balanced separator of size $d(1 + \Delta + \hdots + \Delta^{t})^{-1}$.  This proves the base case. 

Suppose we have a sequence $\S_1, \hdots, \S_i$ with $\S_i \subseteq \S_i'$ and $\beta_i = \bigcap_{S \in \S_1 \cup \hdots \cup \S_i} B(S) \cup C(S)$, such that $\beta_i$ is connected and there exists a normal weight function $w_i$ on $\beta_i$ such that $\beta_i$ has no $(w_i, c)$-balanced separator of size at most $d(1 + \Delta + \hdots + \Delta^{t})^{-i}$. Let $S$ be a separation of $G$ and let $H$ be an induced subgraph of $G$. We define $S \cap H$ as the separation of $H$ given by $(A(S) \cap H, C(S) \cap H, B(S) \cap H)$. If $\S$ is a sequence of separations, we define $\S \cap H = \{S \cap H \text{ s.t. } S \in \S\}$. Let $\S_{i+1}'' = \{S \in S_{i+1}' \text{ s.t. } \cent(S) \subseteq \beta_i\}$. 
Let $T_{i+1} = \{B(S) \cup C(S) : S \in \S_{i+1}''\}$, and let $\S_{i+1}$ be the sequence formed by adding, for every inclusion-wise minimal $Y \in T_{i+1}$, a separation $S \in \S_{i+1}''$ such that $B(S) \cup C(S) = Y$.
Note that since the anchors for canonical separations are contained in their centers, it follows that restricting the anchor map for $\S_{i+1}'$ to $\beta_i$ is a valid anchor map for $\S_{i+1}''$. Note also that by construction, $\S_{i+1}$ is primordial. By Lemma \ref{lemma:primordial-to-Alaminar}, $\S_{i+1}$ is $A$-laminar. Let $\beta_{i+1}$ be the central bag for $\S_{i+1} \cap \beta_i$, and let $w_{i+1}$ be the weight function on $\beta_{i+1}$. Then, 
\begin{align*}
    \beta_{i+1} &=\bigcap_{S \in \S_{i+1} \cap \beta_i} B(S) \cup C(S) \\
    &= \bigcap_{S \in \S_{i+1}} (B(S) \cup C(S)) \cap \beta_i \\
    &= \bigcap_{S \in \S_1 \cup \hdots \cup \S_{i+1}} B(S) \cup C(S). 
\end{align*}

From the choice of $S_{i+1}''$, it holds that $C(S)$ is connected for all $S \in \S_{i+1}''$. Note that since $\S$ is $(a, t)$-good and $d \geq (1 + \Delta + \hdots + \Delta^{t})^k$, it follows that $|C(S)| \leq 1 + \Delta + \hdots + \Delta^t \leq \frac{d}{(1 + \Delta + \hdots + \Delta^{t})^i}$. By Lemmas \ref{lemma:central_bag_1} and \ref{lemma:central_bag_2}, it follows that $\beta_{i+1}$ is connected and that $\beta_{i+1}$ has no $(w_{i+1}, c)$-balanced separator of size $d(1 + \Delta + \hdots + \Delta^{t})^{-(i+1)}$. This completes the induction. Let $\beta = \beta_{k}$, let $w_\S = w_k$, and let $\S^* = \S_1 \cup \hdots \cup \S_k$. Then, $\beta = \bigcap_{S \in \S^*} B(S) \cup C(S)$, $\beta$ is connected, and $w_\S$ is a normal weight function on $\beta$ such that $\beta$ has no $(w_\S, c)$-balanced separator of size at most $d(1 + \Delta + \hdots + \Delta^{t+1})^{-k}$. This proves \eqref{beta-connected} and \eqref{beta-no-balanced-separator}. 

By construction, $\S^* \subseteq \S$, which proves \eqref{S-star-in-S}. It remains to prove (ii). Let $S \in \S \setminus \S^*$, let $\beta_0 = G$ and let $\S_1'' = \S_1'$,  and assume $S \in \S_i'$ for some $1 \leq i \leq k$. Suppose $S \not \in \S_i''$. Then, $\cent(S) \not \subseteq \beta_{i-1}$, so it follows that there exists $S' \in \S_1 \cup \hdots \cup \S_{i-1}$ such that $\cent(S) \cap A(S') \neq \emptyset$. Now, assume $S \in \S_i'' \setminus \S_i$. Because there exists $Y \in T_i$ with $Y \subseteq B(S) \cup C(S)$, it follows that there exists $S' \in \S_{i}$ such that $S'$ is a shield for $S$. This proves \eqref{S-star-is-covering}. 
\end{proof}


    
   
    


Previously, central bags were defined for $A$-laminar sequences of separations. Here, we define central bags for sequences of separations of bounded dimension: we call $\beta$ as in Lemma \ref{lemma:central_bag_high_dimension} a {\em central bag for $\S$}, $w_\S$ the {\em weight function on $\beta$}, and $\S_1, \hdots, \S_k$ the {\em central bag generator for $\S$}. Next, we show how to construct useful sequences of separations of bounded dimension. A sequence $\S$ of separations is {\em strongly laminar} if $C(S_1) \cap C(S_2) = \emptyset$ for all distinct $S_1, S_2 \in \S$. The following lemma states that
under certain conditions, a strongly laminar sequence is laminar.


\begin{lemma}\label{lemma:strongly_laminar}
Let $c \in [\frac{1}{2}, 1)$ and let $\Delta, d$ be positive integers with $d \geq 1 + \Delta + \hdots + \Delta^{t+1}$. Let $G$ be a connected graph with maximum degree $\Delta$, let $w$ be a normal weight function on $G$, and suppose $G$ has no $(w, c)$-balanced separator of size at most $d$. Let $\S$ be a strongly laminar sequence of canonical separations such that $C(S)$ is connected and $|C(S)| \leq d$ for every $S \in \S$. Then, $\S$ is laminar. 
\end{lemma}
\begin{proof}
Let $S_1, S_2 \in \S$. 
Since $C(S_2)$ is connected, it follows that either $C(S_2) \subseteq B(S_1)$ or $C(S_2) \subseteq A(S_1)$. Similarly, either $C(S_1) \subseteq B(S_2)$ or $C(S_1) \subseteq A(S_2)$. Suppose that $C(S_2) \subseteq B(S_1)$ and $C(S_1) \subseteq A(S_2)$. Then, $C(S_2) \cap A(S_1) = \emptyset$ and $C(S_1) \cap B(S_2) = \emptyset$. Since $G$ is connected, it follows that $B(S_2) \cap A(S_1) = \emptyset$, and thus $S_1$ and $S_2$ are non-crossing. The other cases follow by symmetry. 
\end{proof}

Let $X$ be a connected graph, and let $\X = \{Y \subseteq V(G) : Y \text{ is an $X$ in $G$}\}$. Let $\S_X = \{S_Y: Y\in \X\}$ be the sequence of canonical separations with centers in $\X$. We call $\S_X$ the {\em $X$-covering sequence for $G$}. The following lemma shows that if $\S_X$ is $(a, t)$-good, then $\S_X$ has bounded dimension. 

\begin{lemma}
\label{lemma:(a,t)-good_bounded_dimension}
Let $c \in [\frac{1}{2}, 1)$ and let $a, d, t, \Delta$ be positive integers with $d \geq 1 + \Delta + \hdots + \Delta^{t+1}$. Let $G$ be a connected graph with maximum degree $\Delta$, let $w$ be a normal weight function on $G$, and suppose $G$ has no $(w, c)$-balanced separator of size at most $d$. Let $X$ be a connected graph and let $\S_X$ be the $X$-covering sequence. Suppose $\S_X$ is  $(a, t)$-good. Then, $\dim(\S_X) \leq a(1 + \Delta + \hdots + \Delta^{2t}) + 1$. 
\end{lemma}
\begin{proof}
Let $H$ be a graph with $V(H) = \{C(S) : S \in \S_X\}$. Two vertices $C_1, C_2 \in V(H)$ are adjacent if $C_1 \cap C_2 \neq \emptyset$ in $G$. If $C_1 \cap C_2 \neq \emptyset$, then $C_1 \cup C_2$ has diameter at most $2t$, so $\anchor_{\S_X}(S_1)$ has distance at most $2t$ from $\anchor_{\S_X}(S_2)$. Because $G$ has maximum degree $\Delta$, there are at most $1 + \Delta + \hdots + \Delta^{2t}$ vertices in $N^{2t}[v]$ for every $v \in V(G)$. Since every vertex $v \in V(G)$ is the anchor for at most $a$ separations in $\S_X$, it follows that for each $S_1 \in \S_X$ there are at most $a(1 + \Delta + \hdots + \Delta^{2t})$ separations $S_2$ such that $C(S_1) \cap C(S_2) \neq \emptyset$. Therefore, $H$ has maximum degree $a(1 +\Delta + \hdots + \Delta^{2t})$, so $\chi(H) \leq a(1 + \Delta + \hdots + \Delta^{2t}) + 1$. 

Let $\gamma = a(1 + \Delta + \hdots + \Delta^{2t}) + 1$. Let $\chi: V(H) \to \{1, \hdots, \gamma\}$ be a coloring of $H$, 
and let $\S_i = \{S : \chi(C(S)) = i\}$ for all $1 \leq i \leq \gamma$. Now, $\S_i$ is strongly laminar for all $1 \leq i \leq \gamma$, and so by Lemma \ref{lemma:strongly_laminar}, $\S_i$ is laminar for all $1 \leq i \leq \gamma$. Therefore, $\S_1, \hdots, \S_{\gamma}$ is a partition of $\S_X$ into $\gamma$ laminar sequences, so $\dim(\S_X) \leq a(1 + \Delta + \hdots + \Delta^{2t}) + 1$.
\end{proof}


Let $G$ and $X$ be graphs, let $\S_X$ be the $X$-covering sequence in $G$, and let $\beta_X$ be a central bag for $\S_X$. For certain graphs $X$, we can restrict the properties of $\beta_X$ in helpful ways. To do this, we use structures called forcers. Let $G$ be a graph, and let $X, Y \subseteq V(G)$ such that $X \cap Y =\emptyset$.  We say that {\em $X$ breaks $Y$} if for every component $D$ of $G \setminus N[X]$ we have that $Y \not \subseteq N[D]$.
A graph $F$ is
an {\em $X$-forcer for $G$} if for every $Y \subseteq V(G)$ such that $Y$ is an $F$ in $G$, there exists $X' \subset Y$ such that $X'$ is an $X$ in $G$ and 
$X'$ breaks $Y \setminus X'$. For a class $\mathcal{C}$ of graphs, we say
that $F$ is an {\em $X$-forcer for $\mathcal{C}$} if $F$ is an $X$-forcer for
every $G \in \mathcal{C}$.

\begin{lemma}
Let $G$, $X$, and $F$ be graphs such that $F$ is an $X$-forcer for $G$. Let $H$ be an $F$ in $G$. Then, there exists $X'$ which is an $X$ in $G$ such that $X' \subset H$ and $A_{X'} \cap H \neq \emptyset$.
\label{lemma:A_sides_forcer}
\end{lemma}
\begin{proof}
Since $H$ is an $F$ in $G$ and $F$ is an $X$-forcer, there exists $X' \subset H$ such that $X'$ is an $X$ in $H$ and $X'$ breaks $H \setminus X'$. Note that $B_{X'}$ is a component of $G \setminus N[X']$ such that $B_{X'} \cup (C_{X'} \setminus X') \subseteq N[B_{X'}]$. Therefore, $H \setminus X' \not \subseteq B_{X'} \cup C_{X'}$, and so $H \cap A_{X'} \neq \emptyset$. 
\end{proof}


\begin{lemma}
\label{lemma:beta_no_forcers}
Let $c \in [\frac{1}{2}, 1)$ and let $\Delta, d, a, t, k$ be positive integers with $d \geq (1 + \Delta + \hdots + \Delta^{t+1})(1 + \Delta + \hdots + \Delta^{t})^k$. Let $G$ be a connected graph with maximum degree $\Delta$, let $w$ be a normal weight function on $G$, and suppose $G$ has no $(w, c)$-balanced separator of size at most $d$. Let $X$ be a connected graph, let $\S_X$ be the $X$-covering sequence for $G$, and assume $\S_X$ is $(a, t)$-good and $\dim(\S_X) = k$. Let $\beta_X$ be a central bag for $\S_X$. Then, if $F$ is an $X$-forcer for $G$, then $\beta_X$ is $F$-free.  
\end{lemma}
\begin{proof}
Suppose for a contradiction that $H \subseteq \beta_X$ is an $F$ in $G$ and $F$ is an $X$-forcer for $G$. By Lemma \ref{lemma:A_sides_forcer}, there exists $X' \subset H$ such that $X'$ is an $X$ in $G$ and $A_{X'} \cap H \neq \emptyset$. Let $\S_1, \hdots, \S_k$ be the central bag generator for $\S_X$ and let $\S^* = \S_1 \cup \hdots \cup \S_k$. Since $\beta_X \subseteq B(S) \cup C(S)$ for all $S \in \S^*$, it follows that $S_{X'} \not \in \S^*$. Then, by (ii) of Lemma \ref{lemma:central_bag_high_dimension}, there exists $S' \in \S^* \S_1 \cup \hdots \cup \S_k$ such that either $S'$ is a shield for $S_{X'}$ or $X' \cap A(S') \neq \emptyset$. If $X' \cap A(S') \neq \emptyset$, then $H \not \subseteq \beta_X$, a contradiction. Therefore, $S'$ is a shield for $S_{X'}$. But now $A(S_{X'}) \subseteq A(S')$, and $\beta_X \subseteq B(S') \cup C(S')$, so $H \not \subseteq \beta_X$, a contradiction. 
\end{proof}

We now summarize what we have proved so far, as follows.
\begin{theorem}
\label{thm:centralbag}
Let $\mathcal{C}$ be a class of graphs with maximum degree at most $\Delta$ closed under taking induced subgraphs.
Let $t,N$ be integers, and let $X$ be a connected graph with $|V(X)|<t$.
Let $\mathcal{F}$ be a set of graphs such that $F$ is an $X$-forcer for $\mathcal{C}$ for every $F \in \mathcal{F}$. Let $\gamma(x) = 1 + \Delta + \Delta^2 + \hdots + \Delta^x$. If $\tw(H) < N$ for every $\mathcal{F}$-free graph $H$ in
$\mathcal{C}$, then $\tw(G) \leq 2N\gamma(t+1)^{\Delta^{t^2}\gamma(2t) + 1}$ for every
$G \in \mathcal{C}$.
\end{theorem}

\begin{proof}
Let $G \in \mathcal{C}$ and suppose that $\tw(G)> 2N\gamma(t)^{\Delta^{t^2}\gamma(2t) + 1}$. 
We may assume that $G$ is connected.
By Lemma~\ref{lemma:harvey-wood-weights} there exists a uniform weight function
$w$ on $G$ such that $G$ has no $(w,\frac{1}{2})$-balanced separator of size at most $N\gamma(t+1)^{\Delta^{t^2}\gamma(2t) + 1}$. 
Let $\S_X$ be the $X$-covering sequence for $G$. 
Since $|V(X)|<t$, $X$ is connected,  and $G$ has maximum degree $\Delta$, it follows that
every vertex of $G$ belongs to at most ${\Delta^{t} \choose t} \leq {\Delta^{t^2}}$ copies of $X$. Also, since $|V(X)| <t$ and $X$ is connected, for every $S \in \S_X$, $C(S)$ has diameter at most $t$. Therefore $\S_X$ is $(\Delta^{t^2},t)$-good. By Lemma~\ref{lemma:(a,t)-good_bounded_dimension} we have that $\dim(\S_X) \leq \Delta^{t^2}\gamma(2t) + 1$. 
Let $\beta_X$ be a central bag for $\S_X$.
Now by Lemma~\ref{lemma:central_bag_high_dimension} there is a normal
weight function $w_X$ on $\beta_X$ such that $\beta_X$ has no $(w_X,\frac{1}{2})$-balanced separator of size at most $N$.
But by Lemma~ \ref{lemma:beta_no_forcers}, $\beta_X$ is $\mathcal{F}$-free, and therefore
$\tw(\beta_X) < N$,  contrary to
Lemma~\ref{lemma:bounded-tw-balanced-separator}.
\end{proof}

Next, we give a useful application of the results of this section. While it is not used in this paper, it is an important tool for future applications of the central bag method. 

A {\em clique cutset} of a connected graph $G$ is a set $C \subseteq V(G)$ such that $C$ is a clique in $G$ and $G \setminus C$ is not connected. Let $K$ be a clique cutset in $G$, so in particular, $K \neq \emptyset$. The {\em canonical separation for $K$}, denoted $S_K = (A_K, C_K, B_K)$, is defined as follows: $B_K$ is the lexicographically minimum largest-weight connected component of $G \setminus K$, $C_K = K$, and $A_K = V(G) \setminus (B_K \cup C_K)$. A separation is called a {\em clique separation} if it is the canonical separation for some clique cutset $K$ of $G$. Let $\C$ be a primordial sequence of clique separations such that for every clique separation $S_K$ of $G$, it holds that $S_K$ has a shield in $\C$. We call $\C$ a {\em clique covering} of $G$. The next lemma states that $\C$ is $A$-loosely laminar. 

\begin{lemma}
\label{lemma:C-Alooselylaminar}
Let $\Delta, d$ be positive integers with $d > \Delta$ and let $c \in [\frac{1}{2}, 1)$. Let $G$ be a connected graph with maximum degree $\Delta$, let $w$ be a normal weight function on $G$, and suppose $G$ has no $(w, c)$-balanced separator of size at most $d$. Let $\C$ be a clique covering of $G$. Then, $\C$ is $A$-loosely laminar. 
\end{lemma}
\begin{proof}
Suppose there is $S_K, S_{K'} \in \C$ such that $S_K$ and $S_{K'}$ are $A$-loosely crossing. We may assume that $A_{K} \cap C_{K'} \neq \emptyset$. Since $C_{K'}$ is a clique and $A_{K}$ is anticomplete to $B_{K}$, it follows that $C_{K'} \cap B_K = \emptyset$. Since $B_K$ is connected and $A_{K'}$ is anticomplete to $B_{K'}$, it follows that $A_{K'} \cap B_K = \emptyset$. Since $G$ has no $(w, c)$-balanced separator of size $\Delta + 1$, it holds that $w(B_K) > \frac{1}{2}$ and $w(B_{K'}) > \frac{1}{2}$, so $B_K \cap B_{K'} \neq \emptyset$. If $C_K \cap A_{K'} = \emptyset$, then $S_K$ is a shield for $S_{K'}$, a contradiction, so $C_K \cap A_{K'} \neq \emptyset$. Since $C_K$ is a clique and $A_{K'}$ is anticomplete to $B_{K'}$, it follows that $C_K \cap B_{K'} = \emptyset$. Since $B_{K'}$ is connected, it holds that $B_{K'} \cap A_K = \emptyset$, so $B_K = B_{K'}$. Now, $C_K \cap C_{K'}$ is a cutset of $G$ separating $B_K = B_{K'}$ from $A_{K} \cup A_{K'}$, and so, since $G$ is connected, $C_K \cap C_{K'} \neq \emptyset$. But $S_{K''} = (A_{K} \cup A_{K'}, C_{K} \cap C_{K'}, B_K=B_{K'})$ is a shield for both $S_K$ and $S_{K'}$, a contradiction.
\end{proof}

Now, we prove two important results about the central bag for $\C$. 
\begin{theorem}
Let $\Delta, d$ be positive integers with $d > \Delta$ and let $c \in [\frac{1}{2}, 1)$. Let $G$ be a connected graph with maximum degree $\Delta$, let $w$ be a normal weight function on $G$, and assume that $G$ has no $(w,c)$-balanced separator of size at most $d$. Let $\C$ be a clique covering of $G$, let $\beta_\C$ be the central bag for $\C$, and let $w_\C$ be the weight function on $\beta_\C$. Then:
\begin{enumerate}[(i)]
\itemsep -0.2em
\item \label{beta-no-balanced-separator-cliques} $\beta_\C$ has no $(w_\C, c)$-balanced separator of size $d(1+\Delta)^{-1}$, and
\item \label{beta-no-clique-cutset} $\beta_\C$ has no clique cutset. 
\end{enumerate}
\end{theorem}
\begin{proof}
Since cliques have diameter one, it follows from Lemmas \ref{lemma:C-Alooselylaminar} and \ref{lemma:central_bag_2} that $\beta_\C$ has no $(w_\C, c)$-balanced separator of size $d(1+\Delta)^{-1}$. This proves \eqref{beta-no-balanced-separator-cliques}. 

Next, we prove \eqref{beta-no-clique-cutset}. Suppose $\beta_\C$ has a clique cutset $K$. Let $D_1$ and $D_2$ be two connected components of $\beta_C \setminus K$. By Lemma \ref{lemma:central_bag_1} \eqref{C-connected}, $N(A)$ is a clique for every connected component $A$ of $G \setminus \beta_\C$, so we deduce that $D_1$ and $D_2$ are in different connected components of $G \setminus K$. Therefore, $\beta_\C$ intersects two connected components of $G \setminus K$, so $\beta_C \cap A_K \neq \emptyset$. Since $\beta_\C \subseteq B_{K'} \cup C_{K'}$ for all $S_{K'} \in \C$, it follows that $S_K \not \in \C$ and $S_K$ does not have a shield in $\C$, a contradiction. This proves \eqref{beta-no-clique-cutset}.
\end{proof}

\section{Treewidth of claw-free graphs} 
\label{sec:tw_of_clawfree}



	In this section, we prove that for every $k$, every claw-free graph with bounded maximum degree and with no induced subgraph isomorphic to the line graph of a $(k\times k)$-wall has bounded treewidth. Our proof relies on a structural description of claw-free graphs due to the second author and Seymour. In particular, the theorem we apply here is a straightforward corollary of the main result of \cite{Claw5}. To state this theorem, we need a couple of definitions from \cite{Claw5}.
	
	Given a graph $H$, a set $F$ of unordered pairs of vertices of $H$ is called a \emph{valid set} for $H$ if every vertex of $H$ belongs to at most one member of $F$. For a graph $H$ and a valid set $F$ of $H$, we say that a graph $G$ is a \emph{thickening} of $(H,F)$ if for every $v \in V(H)$ there is a nonempty subset $X_v \subseteq V(G)$, all pairwise disjoint and with union $V(G)$, for which the following hold.
	\begin{itemize}
		\itemsep0em
		\item For each $v \in V(H)$, the set $X_v$ is a clique of $G$,
		\item if $u, v \in V(H)$ are adjacent in $H$ and $\{u,v\} \notin F$, then $X_u$ is complete to $X_v$ in $G$,
		\item if $u, v \in V(H)$ are non-adjacent in $H$ and $\{u,v\} \notin F$, then $X_u$ is anticomplete to $X_v$ in $G$,
		\item if $\{u, v\} \in F$, then $X_u$ is neither complete nor anticomplete to $X_v$ in $G$.
	\end{itemize}
	
	Let $\Sigma$ be a circle and let $I = \{I_1, \dots, I_k\}$ be a collection of subsets of $\Sigma$, such that each $I_i$ is homeomorphic to the interval $[0,1]$, no two of $I_1, \dots, I_k$ share an endpoint, and no three of them have union $\Sigma$. Let $H$ be a graph whose vertex set is a finite subset of $\Sigma$, and distinct vertices $u, v \in V(H)$ are adjacent precisely if $u, v \in I_i$ for some $i=1,\ldots, k$. The graph $H$ is called a \emph{long circular interval graph}. Let $F'$ be the set of all pairs $\{u, v\}$ such that $u, v \in V(H)$ are distinct endpoints of $I_i$ for some $i$ and there exists no $j \neq i$ for which $u, v \in I_j$. Also, let $F \subseteq F'$. Then, for every such $H$ and $F$, every thickening $G$ of $(H, F)$ is called a \emph{fuzzy long circular interval graph}.\\

	Given a graph $G$, a \textit{strip-structure} of $G$ is a pair $(H,\eta)$, where $H$ is a graph with no isolated vertices and possibly with loops or parallel edges, and $\eta$ is a function mapping each $e\in E(H)$ to a subset $\eta(e)$ of $V(G)$, and each pair $(e,u)$ consisting of an edge $e\in E(H)$ and an end $u$ of $e$ to a subset $\eta(e,u)$ of $\eta(e)$, with the following specifications.
	\begin{itemize}
		\item[(S1)] The sets $(\eta(e): e\in E(H))$ are non-empty and partition $V (G)$.
		\item[(S2)] For each $v\in V (H)$, the union of sets $\eta(e,v)$ for all $e\in E(H)$ incident with $v$ is a clique of $G$. In particular, $\eta(e,v)$ is a clique of $G$ for all $e\in E(H)$ and $v\in V(H)$ an end of $e$.
		\item[(S3)] For all distinct $e_1,e_2 \in E(H)$, if $x_1\in \eta(e_1)$ and $x_2\in \eta(e_2)$ are adjacent, then there exists $v\in V(H)$ with $v$ an end of both $e_1$ and $e_2$, such that $x_i \in \eta(e_i,v)$ for $i=1,2$.
	\end{itemize}
	We say a strip-structure $(H,\eta)$ is \textit{non-trivial} if $|E(H)|\geq 2$. The following can be derived from Theorem 7.2 in \cite{Claw5}.
	\begin{theorem}[Corollary of Theorem 7.2 from \cite{Claw5}]\label{clawstructure2}
		Let $G$ be a connected claw-free graph. Then one of the following holds.
\begin{itemize}
	\item We have $\alpha(G)\leq 3$.
	\item $G$ is a fuzzy long circular interval graph.
	\item  $G$ admits a non-trivial strip structure $(H,\eta)$, such that for every $e \in E(G)$ with ends $u$ and $v$,
	\begin{itemize}
		\item either  $\alpha( \eta(e))\leq 4$ or $ \eta(e)$ is a fuzzy long circular interval graph; and
		\item there exists a path $P_e$ in $ \eta(e)$ (possibly of length zero) with an end in $\eta(e,u)$ and an end in $\eta(e,v)$ whose interior is disjoint from $\eta(e,u)\cup \eta(e,v)$.
	\end{itemize}
	\end{itemize}
	\end{theorem}
To begin with, we show that every fuzzy long circular interval graph with bounded maximum degree has bounded treewidth. Indeed, the proof is almost immediate from the following well-known fact about \textit{chordal} graphs, i.e. graphs with no induced cycle of length at least four.
\begin{theorem}[folklore]	\label{twchordal}
    	A graph $G$ is chordal if and only if it admits a tree decomposition $(T,\beta)$ where for every $t\in V(T)$, the set $\beta(t)$ is a clique of $G$. Consequently, if $G$ is chordal, then $\tw(G)= \omega(G)-1$. 
\end{theorem}

\begin{theorem}	\label{circtw}
    	Let $G$ be a fuzzy long circular interval graph of maximum degree at most $\Delta$. Then we have $\tw(G) \leq 4\Delta+3$.
\end{theorem}
\begin{proof}
Suppose that $G$ is a thickening of $(H,F)$, where $H$ is a long circular interval graph with $\Sigma, I_1, \dots, I_k$ as in the definition, and $F$ is a valid set for $H$ as in the definition. Let $G^*$ be the graph with $V(G^*)=V(G)$ and 
\[E(G^*)=E(G)\cup \left(\bigcup_{\{u,v\}\in F}\{ab:a\in X_u,b\in X_v\}\right)\cdot\]
Then $G^*$ is a long circular interval graph (the same interval representation $\Sigma, I_1, \dots, I_k$ works for $G^*$, as well). In addition, we may easily observe that
\begin{itemize}
    \item $\omega(G^*)\leq 2\omega(G)\leq 2(\Delta+1)$;  
    \item for all $i=1,\ldots, k$, the set $C_i=\bigcup_{u\in V(H)\cap I_i}X_u$ is a clique of $G^*$; and
    \item for all $i=1,\ldots, k$, the graph $G-C_i$ is a chordal.
\end{itemize}
By Theorem \ref{twchordal} and the third bullet above, $G-C_1$ admits a tree decomposition $(T,\beta)$ of width $\omega(G^*)-1$. Now, for every $t\in V(T)$, let $\beta^*(t)=\beta(t)\cup C_1$. Then it is readily seen that $(T,\beta^*)$ is a tree decomposition of $G^*$ of width $\omega(G^*)+|C_1|-1\leq 2\omega(G^*)-1$, where the last inequality follows from the the second bullet above. Hence, since $G$ is a subgraph of $G^*$, we have $\tw(G)\leq \tw(G^*)\leq 2\omega(G^*)-1\leq 4\Delta+3$, where the last inequality follows from the first bullet above. This proves Theorem \ref{circtw}.
\end{proof}
The following is an easy observation.
\begin{observation}\label{subdtw}
	Let $H$ be a graph and $H'$ be a subdivision of $H$. Then $\tw(H)=\tw(H')$.
\end{observation}
We also use Theorem \ref{wallminor} with an explicit value of $f(k)$. In fact, a considerable amount of work has been devoted to understanding the order of magnitude of $f(k)$, and as of now, the following result of Chuzhoy and Tan provides the best known bound. 
\begin{theorem}[\cite{chuzhoy}]\label{RSwalltw}
	There exist universal constants $c_1$ and $c_2$ such that for every integer $k$, every graph with no subgraph isomorphic to a subdivision of the $(k\times k)$-wall has treewidth at most $c_1k^9\log^{c_2}k$. 
\end{theorem}
Now we are in a position to prove the main result of this section.
\begin{theorem}\label{clawfreelinewalltw}
	Let $\Delta,k$ be integers and $c_1$ and $c_2$ be as in Theorem \ref{RSwalltw}. Let 
	\[w(\Delta,k)=\max\{c_1k^9\log^{c_2}k(\Delta +1)^2,6(\Delta+1)\}-1.\]
	Then for every claw-free graph $G$ of maximum degree at most $\Delta$ and with no induced subgraph isomorphic to the line graph of a subdivision of the $(k\times k)$-wall, we have $\tw(G)\leq w(\Delta,k)$.
\end{theorem}
\begin{proof}
We may assume that $G$ is connected, and so we may apply Theorem \ref{clawstructure2}. Note that if we allow for trivial strip-structures, then the first two bullets of Theorem \ref{clawstructure2} will be absorbed into the first dash of the third bullet. In other words, we have

\sta{\label{pawelclaw}$G$ admits a (possibly trivial) strip structure $(H,\eta)$, such that for every $e \in E(G)$ with ends $u$ and $v$,
	\begin{itemize}
		\item either  $\alpha( \eta(e))\leq 4$ or $ \eta(e)$ is a fuzzy long circular interval graph; and
		\item if $(H,\eta)$ is non-trivial, then there exists a path $P_e$ in $ \eta(e)$ (possibly of length zero) with an end in $\eta(e,u)$ and an end in $\eta(e,v)$ whose interior is disjoint from $\eta(e,u)\cup \eta(e,v)$.
	\end{itemize}}
We also deduce:

\sta{\label{etaclique} For every $e\in E(H)$ and every $v\in V(H)$ incident with $e$, we have $|\eta(e,v)|\leq \Delta+1$.}

By (S2), $\eta(e,v)$ is a clique of $G$. So from $G$ being of maximum degree at most $\Delta$, we have $|\eta(e,v)|\leq \Delta+1$. This proves \eqref{etaclique}.

\sta{\label{Hbnddeg} For every $v\in V(H)$, the number of edges $e\in E(H)$ incident with $v$ for which $\eta(e,v)\neq \emptyset$ is at most $\Delta+1$.}

For otherwise by (S2), the union of sets $\eta(e,v)$ for all $e\in E(H)$ with $v\in e$ contains a clique of $G$ of size at least $\Delta+2$, which is impossible. This proves \eqref{Hbnddeg}.

\sta{\label{Hbnddtw} $H$ admits a tree decomposition $(T_0,\beta_0)$ of width at most $c_1k^9\log^{c_2}k$.}

If $|E(H)|\leq 1$, then we are done. So we may assume that $(H,\eta)$ is non-trivial. Let the paths $\{P_e : e\in E(H)\}$ be as promised in the second bullet of \eqref{pawelclaw}, and let $H^-$ be the graph obtained from $H$ by removing its loops. Then, using (S1), (S2) and (S3) from the definition of a strip structure, one may observe that $G'=G[\bigcup_{e\in E(H^-)}V(P_e)]$ is isomorphic to the line graph of a subdivision $H'$ of $H^-$. Now, since $G$ has no induced subgraph isomorphic to the line graph of a subdivision of the $(k\times k)$-wall, neither does $G'$, and so $H'$ has no subgraph isomorphic to a subdivision of the $(k\times k)$-wall. Thus, by Theorem \ref{RSwalltw}, we have $\tw(H')\leq c_1k^9\log^{c_2}k$, and so by Observation \ref{subdtw}, we have $\tw(H^-)\leq c_1k^9\log^{c_2}k$. This, along with the fact that every tree decomposition of $H^-$ is also a tree decomposition of $H$, proves \eqref{Hbnddtw}.

\sta{\label{etabnddtw} For every $e\in E(H)$, $ \eta(e)$ admits a tree decomposition $(T_e,\beta_e)$ of width at most $4(\Delta+1)$.}

Note that $ \eta(e)$ is of maximum degree at most $\Delta$. So if $\alpha( \eta(e))\leq 4$, then we have $\tw( \eta(e))\leq |\eta(e)|\leq \alpha( \eta(e))(\Delta+1)\leq 4(\Delta+1)$, as desired. Otherwise, by the first bullet of \eqref{pawelclaw}, $ \eta(e)$ is a fuzzy long circular interval graph, and so by Theorem \ref{circtw}, we have $\tw(G)\leq 4\Delta+3$. This proves \eqref{etabnddtw}.\vsp
 
Let $(T_0,\beta_0)$ be as in \eqref{Hbnddtw}, and for every $e\in E(H)$, let $(T_e,\beta_e)$ be as promised by \eqref{etabnddtw}. We assume $T_0$ and $T_e$'s have mutually disjoint vertex sets and edge sets. Now, we construct a tree $T$ as follows. For every $e \in E(H)$ with ends $u$ and $v$, choose a vertex $s_e\in \beta^{-1}_0(u)\cap \beta^{-1}_0(v)$, which exists by definition of tree decomposition, and pick $t_e\in V(T_e)$ arbitrarily. Let $V(T)=V(T_0)\cup (\bigcup_{e\in E(H)}V(T_e))$, and $E(T)= \{s_et_e:e\in E(H)\}\cup E(T_0)\cup (\bigcup_{e\in E(H)}E(T_e))$. We also define $\beta:V(T)\rightarrow 2^{V(G)}$ as follows. Let $t\in V(T)$. If $t\in V(T_0)$, then
 \[\beta(t)=\bigcup_{u \in \beta_0(t)} \bigcup_{e \in E(H):\text{$u$ is an end of $e$}} \eta(e,u) \cdot\] Otherwise, if $t\in V(T_e)$ for some $e \in E(H)$ with ends $u$ and $v$, then $\beta(t)=\beta_e(t)\cup \eta(e,u) \cup \eta(e,v)$.
 
 \sta{\label{finaltd} $(T,\beta)$ is a tree decomposition of $G$.}
By (S1), for every vertex $x\in V(G)$, there exists $e\in E(H)$ such that $x\in \eta(e)$, and so $(T_e,\beta_e)$ being a tree decomposition of $ \eta(e)$, there exists $t\in V(T_e)\subseteq V(T)$ with $x\in \beta_{e}(t)\subseteq \beta(t)$.

Also, by (S3), for every edge $x_1x_2\in E(G)$, either $x_1x_2\in E( \eta(e))$ for some $e\in E(H)$, or there exists $v\in V(H)$ and $e_1,e_2\in E(H)$ with $v$ an end of $e_1$ and $e_2$ such that $x_i\in \eta(e_i,v)$ for $i=1,2$. In the former case, since $(T_e,\beta_e)$ is a tree decomposition of $ \eta(e)$, there exists $t\in V(T_e)\subseteq V(T)$ with $x_1,x_2\in \beta_{e}(t)\subseteq \beta(t)$. In the latter case, since $(T_0,\beta_0)$ is a tree decomposition of $H$, there exists $t\in V(T_0)\subseteq V(T)$ with $v\in \beta_0(t)$. Therefore, for  $i=1,2$, we have
 \[x_i\in \eta(e_i,v)\subseteq \bigcup_{e \in E(H):\text{ } v \text{ is an end of $e$}} \eta(e,v)\subseteq \bigcup_{u \in \beta_0(t)} \bigcup_{e \in E(H):\text{ } u \text{ is an end of $e$}} \eta(e,u)=\beta(t),\]
and so $x_1,x_2\in  \beta(t)$. 

It remains to show that for every $x\in V(G)$, the graph $T|\beta^{-1}(x)$ is connected. By (S1), there exists a unique edge $e \in E(H)$ with ends $u$ and $v$, with $x\in \eta(e)$. First, suppose that either $x\in \eta(e,u)$ or $x\in \eta(e,v)$, say the former. Then we have $\beta^{-1}(x)=\beta_0^{-1}(u)\cup V(T_e)$. Also, since $s_e\in \beta_0^{-1}(u)$ and $t_e\in V(T_e)$, we have $E(T|\beta^{-1}(x))=\{s_et_e\}\cup E(T_0|\beta_0^{-1}(u))\cup E(T_e)$. Now, from $(T_0,\beta_0)$ being a tree decomposition of $H$,  we deduce that $T_0|\beta_0^{-1}(u)$ is connected, and so $T|\beta^{-1}(x)$ is connected, as well.

Next, suppose that $x\in \eta(e)\setminus (\eta(e,u)\cup  \eta(e,v))$. Then we have $\beta^{-1}(x)=\beta_e^{-1}(x)$. So from $(T_e,\beta_e)$ being a tree decomposition of $ \eta(e)$, we deduce that $T|\beta^{-1}(x)=T_e|\beta_e^{-1}(x)$ is connected. This proves \eqref{finaltd}.\vsp

Now, let $t\in V(T)$. If $t\in V(T_0)$, then by \eqref{etaclique}, \eqref{Hbnddeg} and \eqref{Hbnddtw}, we have
\[|\beta(t)|=\sum_{u \in \beta_0(t)} \sum_{e \in E(H):\text{ $u$ is an end of $e$}} |\eta(e,u)|\leq |\beta_0(t)|(\Delta+1)^2\leq c_1k^9\log^{c_2}k(\Delta +1)^2\leq w(\Delta,k)+1\cdot\]
Also, if $t\in V(T_e)$ for some $e\in E(H)$, then by \eqref{etaclique} and \eqref{etabnddtw}, we have
\[|\beta(t)|\leq |\beta_e(t)|+|\eta(e,u)|+|\eta(e,v)|\leq 6(\Delta+1)\leq w(\Delta,k)+1\cdot\]
Hence, by \eqref{finaltd}, $(T,\beta)$ is a tree decomposition of $G$ of width at most $w(\Delta,k)$. This proves Theorem \ref{clawfreelinewalltw}.
\end{proof}
	

\section{Long claws and line graphs of walls}
\label{sec:claw_free_result}

Here, we apply the  results of Sections \ref{sec:tw_of_clawfree} to prove Theorem \ref{thm:claw-free_nonspecific}, that excluding a long claw and the line graphs of all subdivisions of $W_{k \times k}$ gives bounded treewidth.

Let $t_1, t_2, t_3$ be integers, with $t_1 \geq 0$ and $t_2, t_3 \geq 1$. Recall from the introduction that a {\em long claw}, also called a {\em subdivided claw}, denoted $S_{t_1, t_2, t_3}$, is a vertex $v$ and three paths $P_1$, $P_2$, $P_3$, of length $t_1$, $t_2$, and $t_3$, respectively, with one end $v$, such that $P_1 \setminus \{v\}$, $P_2 \setminus \{v\}$, and $P_3 \setminus \{v\}$ are pairwise disjoint and anticomplete to each other. We call $P_1, P_2, P_3$ the {\em paths} of $S_{t_1, t_2, t_3}$. The vertex $v$ is called the {\em root} of $S_{t_1, t_2, t_3}$.
For two graphs  $H_1,H_2$, we denote by $H_1+H_2$ the graph with vertex set
$V(H_1) \cup V(H_2)$ and edge set $E(H_1) \cup E(H_2)$.
We start with a lemma.

\begin{lemma}
Let $t_1, t_2, t_3$ be positive integers with $t_1 \geq 2$. Let $G$ be an $S_{t_1, t_2, t_3}$-free graph. Then, $S_{t_1-1, t_2, t_3} + K_1$ is an $S_{t_1-2, t_2, t_3}$-forcer for $G$. 
\label{lemma:claw_forcers}
\end{lemma}
\begin{proof}
Let $H$ be an $S_{t_1-1, t_2, t_3}$ in $G$, and let $u \in V(G)$ be anticomplete to $H$, so that $H \cup \{u\}$ is an $S_{t_1-1, t_2, t_3} + K_1$. Let $H = P_1 \cup P_2 \cup P_3$, where $P_1 = v \dd x_1 \dd \hdots \dd x_{t_1-1}$, $P_2 = v \dd y_1 \dd \hdots \dd y_{t_2}$, and $P_3 = v \dd z_1 \dd \hdots \dd z_{t_3}$. 
Let $X = H \setminus x_{t_1-1}$. Let $D$ be a connected component of $G \setminus N[X]$. Suppose $u, x_{t_1 -1} \in N[D]$. It follows that there exists a path $P = x_{t_1 - 1} \dd p_1 \dd \hdots \dd p_k \dd u$ from $x_{t_1 - 1}$ to $u$ with $P^* \subseteq D$, so $X$ is anticomplete to $P^*$. Then, $H \cup \{p_1\}$ is isomorphic to $S_{t_1, t_2, t_3}$, a contradiction.  Therefore, $X$ breaks $\{u, x_{t_1 - 1}\}$, and it follows that $S_{t_1-1, t_2, t_3} + K_1$ is an $S_{t_1-2, t_2, t_3}$-forcer for $G$. \end{proof}

Now we can prove Theorem \ref{thm:claw-free_nonspecific}, which we restate. 

\begin{theorem}
\label{thm:claw_free}
Let $\Delta, t_1,t_2,t_3, k$ be positive integers with $t = t_1+t_2+t_3$. Let  $\mathcal{C}$ be the class of all $S_{t_1, t_2, t_3}$-free graph with maximum degree $\Delta$ and no induced subgraph isomorphic to the line graph of a subdivision of $W_{k \times k}$. There exists an integer $N_{k,t,\Delta}$ such that $\tw(G) \leq N_{k,t,\Delta}$ for every $G \in \mathcal{C}$.
\end{theorem}

\begin{proof}
  The proof is by induction on $t_1+t_2+t_3$. If $t_1=t_2=t_3=1$, the result
  follows from Theorem~\ref{clawfreelinewalltw}. Thus we may assume that
  $t_1 \geq 2$.
  By Theorem~\ref{thm:centralbag} and Lemma~\ref{lemma:claw_forcers}, it is enough to find a bound on the treewidth of $(S_{t_1-1, t_2, t_3}+K_1)$-free graphs in $\mathcal{C}$.

  Let $H \in \mathcal{C}$ be  $(S_{t_1-1, t_2, t_3}+K_1)$-free.
  By the inductive hypothesis we may assume that there exists
  $X \subseteq V(H)$ such that $X$ is an 
  $S_{t_1-1, t_2, t_3}$ in $H$.
  Since $H$ does not contain $S_{t_1-1, t_2, t_3} + K_1$, it follows that
  $V(H) \subseteq N[X]$, and therefore $\tw(H) \leq |V(H)| \leq t\Delta$.
  \end{proof}

\section{$t$-thetas, $t$-pyramids, and line graphs of walls}
\label{sec:theta_pyramid}
In this section, we prove Theorem \ref{thm:pyramid_theta-nonspecific}, that for all $k,t$, excluding $t$-thetas, $t$-pyramids, and the line graphs of all subdivisions of $W_{k \times k}$ in graphs with bounded degree gives bounded treewidth. The proof involves an application of Theorem \ref{thm:claw_free}. We also need the following lemma. 

\begin{lemma}
Let $x_1, x_2, x_3$ be three distinct vertices of a graph $G$. Assume that $H$ is a connected induced subgraph of $G \setminus \{x_1, x_2, x_3\}$ such that $H$ contains at least one neighbor of each of $x_1$, $x_2$, $x_3$, and that subject to these conditions $V(H)$ is minimal subject to inclusion. Then, one of the following holds:
\begin{enumerate}[(i)]
\item For some distinct $i,j,k \in  \{1,2,3\}$, there exists $P$ that is either a path from $x_i$ to $x_j$ or a hole containing the edge $x_ix_j$ such that
\begin{itemize}
\item $H = P \setminus \{x_i,x_j\}$, and
\item either $x_k$ has at least two non-adjacent neighbors in $H$ or $x_k$ has exactly two neighbors in $H$ and its neighbors in $H$ are adjacent.
\end{itemize}

\item There exists a vertex $a \in H$ and three paths $P_1, P_2, P_3$, where $P_i$ is from $a$ to $x_i$, such that 
\begin{itemize}
\item $H = (P_1 \cup P_2 \cup P_3) \setminus \{x_1, x_2, x_3\}$, and 
\item the sets $P_1 \setminus \{a\}$, $P_2 \setminus \{a\}$ and $P_3 \setminus \{a\}$ are pairwise disjoint, and
\item for distinct $i,j \in \{1,2,3\}$, there are no edges between $P_i \setminus \{a\}$ and $P_j \setminus \{a\}$, except possibly $x_ix_j$.
\end{itemize}

\item There exists a triangle $a_1a_2a_3$ in $H$ and three paths $P_1, P_2, P_3$, where $P_i$ is from $a_i$ to $x_i$, such that
\begin{itemize}
\item $H = (P_1 \cup P_2 \cup P_3) \setminus \{x_1, x_2, x_3\} $, and 
\item the sets $P_1$, $P_2$ and $P_3$ are pairwise disjoint, and
\item for distinct $i,j \in \{1,2,3\}$, there are no edges between $P_i$ and $P_j$, except $a_ia_j$ and possibly $x_ix_j$.
\end{itemize}
\end{enumerate}
\label{lem:three_leaves}
\end{lemma}
\begin{proof}
 For some distinct $i,j,k \in  \{1,2,3\}$, let $P$ be a path from $x_i$ to $x_j$ with $V(P^*) \subseteq V(H)$ (in the graph where the edge $x_ix_j$ is deleted if it exists). Such a path exists since $x_i$ and $x_j$ have neighbors in $H$ and $H$ is connected. Assume that $x_k$ has neighbors in $P^*$. Then, by the minimality of $V(H)$, we have $H = P^*$. If $x_k$ has two non-adjacent neighbors in $P^*$, or $x_k$ has two neighbors in $P^*$ and its neighbors in $P^*$ are adjacent, then outcome (i) holds. If $x_k$ has a unique neighbor in $P^*$, then outcome (ii) holds. Thus, we may assume that $x_k$ is anticomplete to $P^*$.

Let $Q$ be a path with $Q \setminus \{x_k\} \subseteq H$ from $x_k$ to a vertex $w \in H \setminus P$ (so $x_k \neq w$) with a neighbor in $P^*$. Such a path exists since $x_k$ has a neighbor in $H$, $x_k$ is anticomplete to $P^*$, and $H$ is connected. By the minimality of $V(H)$, we have $V(H) = (V(P) \cup V(Q)) \setminus \{x_1, x_2, x_3\}$ and no vertex of $Q \setminus w$ has a neighbor in $P^*$. Moreover, by the argument of the previous paragraph, we may assume  that $x_i$ and $x_j$ are anticomplete to $Q \setminus \{x_k\}$.

Now, if $w$ has a unique neighbor in $P^*$, then outcome (ii) holds. If $w$ has two neighbors in $P^*$ and its neighbors in $P^*$ are adjacent, then outcome (iii) holds. Therefore, we may assume that $w$ has two non-adjacent neighbors in $P^*$. Let $y_i$ and $y_j$ be the neighbors of $w$ in $P^*$ that are closest in $P^*$ to $x_i$ and $x_j$, respectively. Let $R$ be the subpath of $P^*$ from $y_i$ to $y_j$. Now, the graph $H'$ induced by $\left((P \cup Q) \setminus R^* \right) \setminus \{x_1, x_2, x_3\}$ is a connected induced subgraph of $G \setminus \{x_1, x_2, x_3\}$ and it contains at least one neighbor of $x_1$, $x_2$, and $x_3$. Moreover, $H'  \subset H$ since $R^* \neq \emptyset$. This contradicts the minimality of $V(H)$.
\end{proof}

Now we are ready to prove Theorem \ref{thm:pyramid_theta-nonspecific}, which we restate.
\begin{theorem}
  Let $\Delta, t,k$ be positive integers with $t \geq 2$. Let $\mathcal{C}$ be the class of 
  graphs of maximum degree $\Delta$ with no $t$-theta, no $t$-pyramid, and no induced subgraph isomorphic to the line graph of a subdivision of $W_{k \times k}$. There exists an integer $M_{k,t,\Delta}$ such that $\tw(G) \leq M_{k,t,\Delta}$ for every $G \in \mathcal{C}$.
\label{thm:theta_pyramid}
\end{theorem}
\begin{proof}
We start by proving a result about the existence of forcers for $\mathcal{C}$.

\sta{\label{claws_are_forcers} $S_{t, t, t}$ is an $S_{t-1, t-1, t-1}$-forcer for
  $\mathcal{C}$.}

Let $G \in \mathcal{C}$, and let $Y$ be an $S_{t, t, t}$ in $G$, let $r$ be the root of $Y$, let $x, y, z$ be the leaves of $Y$, and let $X = Y \setminus \{x, y, z\}$. Let $D$ be a connected component of $G \setminus N[X]$, and suppose $\{x, y, z\} \subseteq N[D]$. Let $Z \subseteq D$ be an inclusion-wise minimal connected subset of $D$ such that $x, y, z$ each have a neighbor in $Z$. By Lemma \ref{lem:three_leaves}, one of three cases holds. If case (ii) or case (iii) holds, then it is clear that $Y \cup Z$ is either a $t$-theta or a $t$-pyramid, so we may assume case (i) holds. Then, up to symmetry between $x, y$, and $z$, the subgraph of $G$ induced on $Z \cup \{x, z\}$ is a path from $x$ to $z$. Suppose $y$ has two non-adjacent neighbors in $Z$. Let $p, q$ in $Z$ be the first and last neighbors of $y$ in $Z$, such that $x, p, q, z$ appear in $x \dd Z \dd z$ in that order. Then $G$ contains a theta between $r$ and $y$ through $r \dd Y \dd y$, $r \dd Y \dd x \dd Z \dd p \dd y$, and $r \dd Y \dd z \dd Z \dd q \dd y$. Since each of the paths of the theta contains a path of $Y$, it follows that every path of the theta has length at least $t$, a contradiction. Therefore, $y$ has exactly two adjacent neighbors $p, q$ in $Z$ such that $x, p, q, z$ appear in $x \dd Z \dd z$ in that order. But now $G$ contains a pyramid from $r$ to $\{y, p, q\}$ through $r \dd Y \dd y$, $r \dd Y \dd x \dd Z \dd p$, and $r \dd Y \dd z \dd Z \dd q$. Since each of the paths of the pyramid contains a path of $Y$, it follows that every path of the pyramid has length at least $t$, a contradiction. Therefore, $X$ breaks $\{x, y, z\}$, so $S_{t, t, t}$ is an $S_{t-1, t-1, t-1}$-forcer for $G$. This proves \eqref{claws_are_forcers}. \vsp

Now by Theorem~\ref{thm:centralbag}, the result follows immediately
from Theorem~\ref{thm:claw_free}.
\end{proof}

\section{Subcubic subdivided caterpillars and their line graphs}\label{sec:caterpillar}
In this section, we prove Theorem \ref{thm:caterpillar-non_specific}, that excluding a subdivided subcubic caterpillar and its line graph in graphs with bounded degree gives bounded treewidth. The proof uses Theorem \ref{thm:claw_free} to get a structural result involving a family of induced subgraphs called $(k, t)$-creatures. We begin with the following lemma. 

\begin{lemma}\label{subpath}
		Let $\Delta>0$ and $t> 0$ be integers, $G$ be a graph and $P$ be an induced path in $G$ of length at least $t(1+\Delta)-1$. Also, let $z\in G\setminus P$ have at least one and at most $\Delta$ neighbors in $P$. Then there exists a subpath $P'=p'_0\dd \cdots\dd p'_{t}$ of $P$ of length $t$ where $N(z)\cap P'=\{p'_0\}$.
	\end{lemma}
	\begin{proof}
		Suppose not. Let $P=p_0\dd\cdots\dd p_{\ell}$, where $\ell\geq t(1+\Delta)-1$. Also, let $|N(z)\cap P|=j\leq \Delta$, and $0\leq i_1<\cdots<i_j\leq \ell$ satisfy $N(z)\cap P=\{p_{i_k}:k=1,\ldots, j\}$. If the subpath $p_{i_1}\dd\cdots\dd p_{0}$ of $P$ is of length at least $t$, then $P'=p_{i_1}\dd \cdots\dd p_{i_1-t}$ satisfies Lemma \ref{subpath}, a contradiction. So $p_{i_1}\dd \cdots\dd p_{0}$ is of length at most $t-1$. Similarly, $p_{i_j}\dd \cdots\dd p_{\ell}$ is of length at most $t-1$. As a result, $j\geq 2$.
		
		Now, if for some $k\in \{1,\ldots, j-1\}$, the subpath $p_{i_k}\dd \cdots\dd p_{i_{k+1}}$ of $P$ is of length at least $t+1$, then $P'= p_{i_k}\dd \cdots\dd p_{i_{k}+t}$ satisfies the lemma. Thus,  for all $k\in \{1,\ldots, j-1\}$, $p_{i_k}\dd \cdots\dd p_{i_{k+1}}$ is of length at most $t$. But then $P$ is of length at most $2(t-1)+t(j-1)=t(j+1)-2\leq t(1+\Delta)-2$, which is impossible. This proves Lemma \ref{subpath}.
	\end{proof}
	Next, we define creatures properly. For integers $k>0$ and $t\geq 0$, a  \emph{$(k,t)$-creature} in a graph $G$ is a pair $\Xi=(J,\mathcal{P})$, where 
	\begin{itemize}
		\item $J$ is a connected subset of $G$.
		\item $\mathcal{P}$ is a collection of $k$ mutually vertex-disjoint and anticomplete induced paths in $G \setminus J$, each of length $t$.
		\item For every $P\in \mathcal{P}$, an end $v$ of $P$, called the \textit{$P$-joint} of $\Xi$, satisfies the following:
		\begin{itemize}
			\item $v$ has a neighbor in $J$, and
			\item $P\setminus v$ is anticomplete to $J$.
		\end{itemize} 
	\end{itemize}
	We also use $\Xi$ to denote the set $J\cup (\bigcup_{P\in \mathcal{P}}P)$.
	\begin{lemma}\label{creature}
		Let $\Delta,k,t\geq 0$ be integers and let $G$ be a graph of maximum degree at most $\Delta$ with no $(k,t)$-creature. Let $X\subseteq G$ be an $S_{t+1,t+1,t+1}$in $G$ and let $x\in X$ be a leaf of $X$. Then the connected component of $G\setminus (N[X\setminus \{x\}]\setminus \{x\})$ containing $x$ has no $(k-1,t(1+\Delta))$-creature. 
		\end{lemma}
	\begin{proof}
		Suppose not. Since $S_{t+1, t+1, t+1}$ is a $(3, t)$-creature, it follows that $k \geq 4$. So we may choose a $(k-1,t(1+\Delta))$-creature $\Xi=(J,\{P_1,\ldots,P_{k-1}\})$ in the component $C$ of $G \setminus (N[X \setminus \{x\}] \setminus \{x\})$ containing $x$ with $P_i$-joint $v_i$ for $i=1,\ldots, {k-1}$, and an induced path $L$ in $C$ from $x$ to some vertex $z\in N[\Xi] $, such that no vertex in $L\setminus \{z\}$ belongs to $N[\Xi]$.
		Let $u$ be the root of $X$ and let $P, Q, R$ be the paths of $X$, with $x$ an end of $P$. 
		We deduce the following.
		
		\sta{\label{1}$z$ has a neighbor in $\bigcup_{i=1}^{k-1} P_i$.}
		
		Suppose for a contradiction that $z$ is anticomplete to $\bigcup_{i=1}^{k-1}P_i$, and so $z\in N[J]$. Note that the path $z\dd L\dd x\dd P\dd u$ has length at least $t$ (since $P$ does), and so we may choose a subpath $P'$ of $z\dd L\dd x\dd P\dd u$ containing $z$ and of length equal to $t$. Also, for each $i\in \{1,\ldots, k-1\}$, since $P_i$ is of length $t(1+\Delta)\geq t$, we may choose a subpath $P'_i$ of $P_i$ containing $v_i$ and of length equal to $t$. But then $(J,\{P'_i: i\in \{1,\ldots, k-1\}\}\cup \{P'\})$ is a $(k,t)$-creature in $G$, a contradiction. This proves \eqref{1}.\vsp
		
		Let $I$ be the set of all indices $i\in \{1,\ldots,k-1\}$ for which $z$ has a neighbor in $P_i$. By \eqref{1}, we have $I\neq \emptyset$, and so we may select an element $i_0\in I$. Let $J'=J\cup P_{i_0}\cup L\cup P$. Note that $J'$ is connected, and $u,z\in J'$. 
		
		\sta{\label{2}For each $i\in I\setminus\{i_0\}$, there exists a subpath $P'_i=p^i_0\dd \cdots\dd p^i_{t}$ of $P_i$ of length $t$, such that $p^i_0$ is the only neighbor of $z$ in $P'_i$, and $P'_i$ is anticomplete to $J'\setminus \{z\}$.}
		
		If $t=0$, then $P_i=\{v_i\}$, $z$ is adjacent to $v_i$, and so $P'_i=P_i$ satisfies \eqref{2}. Thus, we may assume that $t\geq 1$, and so $P_i$ is of length at least one. By Lemma \ref{subpath} applied to $P_i\setminus v_i$ and $z$, we obtain a subpath $P'_i=p^i_0\dd \cdots\dd p^i_{t}$ of $P_i\setminus v_i$ (and so of $P_i$) of length $t$, such that $p^i_0$ is the only neighbor of $z$ in $P'_i$. Also, by the choice of $X$, $\Xi$ and $L$, we obtain that $P_i\setminus \{v_i\}$ is anticomplete to $J'\setminus \{z\}$. Therefore, $P'_i\subseteq P_i\setminus \{v_i\}$ is anticomplete to $J'\setminus \{z\}$, as well. This proves \eqref{2}. \vsp
		
		Now, for each $i\in I\setminus\{i_0\}$, let $P'_i$ be as promised in \eqref{2}. Moreover, for each $i\in \{1,\ldots, k-1\}\setminus (I\setminus \{i_0\})$, there exists a subpath $P'_i$ of $P_i$ containing $v_i$ and of length equal to $t$, as $P_i$ is of length $t(1+\Delta)\geq t$. But then by \eqref{2} and the choice of $X$, $\Xi$ and $L$, $(J',\{P'_i:i\in \{1,\ldots, k-1\}\setminus \{i_0\}\}\cup \{Q\setminus u,R\setminus u \})$ is a $(k,t)$-creature in $G$, a contradiction. This proves Lemma \ref{creature}.
	\end{proof}

\begin{lemma}\label{ramseypath}
Let $\Delta>0$ and $\ell>1$ be integers, and $K$ be a connected graph of maximum degree at most $\Delta$ with $|K|\geq 1+\sum_{i=0}^{\ell-2}\Delta^i$. Then $K$ contains an induced path on at least $\ell$ vertices. 
\end{lemma}
\begin{proof}
    Since $G$ has maximum degree at most $\Delta$, then for every $i\geq 0$ and every vertex $v\in G$, the set of vertices in $G$ at distance $i$ from $v$ is of size at most $\Delta^{i}$. Therefore, $G$ has a vertex at distance at least $\ell-1$ from $v$. This proves Lemma \ref{ramseypath}.
\end{proof}

Recall from the introduction that a tree $T$ is a {\em subcubic subdivided caterpillar} if it is of maximum degree at most three, and there exists a path $P \subseteq T$ such that $P$ contains every vertex of $T$ of degree three. The {\em spine} of $T$ is the shortest path containing all vertices of degree at least three in $T$. A {\em leg} of a subdivided caterpillar $T$ is a path in $T$ from a leaf to a vertex of degree three in $T$ whose all internal vertices are of degree two.

	\begin{theorem}\label{cater}
		Let $T$ be a subcubic subdivided caterpillar and let $\Delta>0$ be an integer. Then there exist $k,t$ such for every graph $G$ of maximum degree at most $\Delta$, if $G$ contains a $(k,t)$-creature, then $G$ contains a subdivision of $T$ or the line graph of a subdivision of $T$.
	\end{theorem}
\begin{proof}
Note that if $T$ has no vertex of degree three, then it is a path, and so setting $k=1$ and $t=|T|$, we are done. So we may assume that $S$ has a spine $S$ with $|S|=s\geq 1$. We define $\ell=6s^3+s^2-1>1$ and $k= 1 + \Delta + \Delta^2 + \hdots + \Delta^{\ell - 2}$. Also, for every leaf $u$ of $T$, let $U_{u}$ be the leg of $T$ having $u$ as one of its end.  Let $t$ be the maximum length of $U_{u}$ taken over all leaves $u \in T\setminus S$ of $T$.
We claim that the values of $k,t$ defined as above satisfy the theorem. Suppose not. Then there exists a graph $G$ of maximum degree at most $\Delta$, containing a $(k,t)$-creature but not containing a subdivision of $T$ or the line graph of a subdivision of $T$.

\sta{\label{12} We may choose $H$ and $\mathcal{P}$ with the following specifications.
\begin{itemize}
\item $H$ is a connected induced subgraph of $G$.
\item $\mathcal{P}$ is a collection of $k$ mutually vertex-disjoint and anticomplete induced paths in $G$, each of length at least $t$. 
\item For every $P\in \mathcal{P}$, there is an end of $P$, denoted by $v_P$, with $P\cap H=\{v_P\}$ and $N(P\setminus \{v_P\})\cap H=\{v_P\}$.
\end{itemize}}
Note that there exists a $(k,t)$-creature $\Xi=(J,\mathcal{P})$ in $G$. For every $P\in \mathcal{P}$, let $v_P$ be the $P$-joint of $\Xi$. Let $H=(J\cup \{v_P:P\in \mathcal{P}\})$. From the definition of a $(k,t)$-creature, it follows directly that $H$ and $\mathcal{P}$ satisfy the above three bullets. This proves \eqref{12}.\vsp

We choose $H$ and $\mathcal{P}$ satisfying \eqref{12} and with $|H|$ as small as possible. For every $P\in \mathcal{P}$, let $v_P$ be as in the third bullet of \eqref{12}. Let $A=\{v_P:P\in \mathcal{P}\}$ and $J=H\setminus A$.

\sta{\label{3} Every vertex in $v\in J$ is a cut-vertex of $H$.}

For otherwise $H\setminus v$ and $\mathcal{P}$ satisfy \eqref{12}, violating the minimality of $H$. This proves \eqref{3}.\vsp

For every vertex $v\in H$, let us say $v$ is \textit{redundant} if $v\in J$ and $N(v)\cap H$ is a stable set in $H$ of size exactly two. Otherwise, we say $v$ is  \textit{irredundant}.

\sta{\label{4}  There exists an induced path $Q_1$ in $H$ containing at least $6s^3+s^2-1$ irredundant vertices.}

For every redundant vertex $v\in H$ with $N(v)\cap H=\{x,z\}$, by \textit{suppressing} $v$, we mean removing $v$ from $H$ and adding the edge $xz$ to the resulting graph, while we refer to the reverse operation as \textit{unsuppressing} $v$. Let $K$ be the graph obtained from $H$ by repeatedly suppressing redundant vertices until none is left. Note that the maximum degree of $K$ does not exceed that of $H$, which in turn does not exceed $\Delta$, as $H$ is an induced subgraph of $G$. Also, we have $A\subseteq K$, and so $|K|\geq |A|=k=1+\sum_{i=0}^{\ell-2}\Delta^i$. Thus, by Lemma~\ref{ramseypath}, $K$ contains an induced path $Q_0$ on at least $\ell=6s^3+s^2-1$ vertices. Therefore, after unsupressing all redundant vertices of $H$, we obtain an induced path $Q_1$ in $H$ where every vertex in $Q_0$ is an irredundant vertex of $Q_1$. This proves \eqref{4}.\vsp

Henceforth, let $Q_1$ be as guaranteed in \eqref{4}.

\sta{\label{5} There exists an induced path $Q_2$ in $H$ such that $Q_2\subseteq J$ and $Q_2$ contains at least $6s^2+s-1$ irredundant vertices.}

Let $A_1=Q_1\cap A$ and $B_1=Q_1\cup (\bigcup_{P\in \mathcal{P}: v_P\in A_1}P)$. If $|A_1|\geq s$, then $B_1$ contains a subdivision of $T$, and hence so does $G$, a contradiction. It follows that $|A_1|\leq  s-1$. As a result, $Q_1\setminus A_1$ has at most $s$ connected components, and by \eqref{4}, $Q_1\setminus A_1$ contains at least $6s^3+s^2-s$ irredundant vertices. Therefore, there exists a connected component $Q_2$ of $Q_1\setminus A_1$ (hence $Q_2\subseteq J$) containing at least $6s^2+s-1$ irredundant vertices. This proves \eqref{5}.\vsp 

Henceforth, let $Q_2$ be as promised in \eqref{5}. Note that by \eqref{1}, every vertex in $Q_2$ is a cut-vertex of $H$. We say a vertex $x\in Q_2$ is \textit{docile} if there exists a connected component of $H\setminus x$, denoted by $D_x$, such that $D_x\cap Q_2=\emptyset$. The following is immediate from the definition.

\sta{\label{6} Let $x\in Q_2$ be a docile vertex. Then 
\begin{itemize}
		\item $D_x$ is anticomplete to $Q_2\setminus \{x\}$;
		\item $N(x)\cap D_x\neq \emptyset$; and
	    \item for every docile vertex $y\in Q_2\setminus \{x\}$,  $D_x$ is anticomplete to $D_y$ in $G$.
	\end{itemize}}
	
Also, we deduce:

\sta{\label{7}  For every docile vertex $x\in G$, there exists a (possibly not unique) path $P_x\in \mathcal{P}$ with $v_{P_x}\in D_x$.}

Otherwise $H\setminus D_x$ and $\mathcal{P}$ satisfy \eqref{12}, violating the minimality of $H$. This proves \eqref{7}.

\sta{\label{8}  There is a subpath $Q_3$ of $Q_2$ which has at least $6s$ irredundant vertices and no docile vertices.}

Let $D$ be the set of all docile vertices in $Q_2$. For every $x\in D$, let $P_x$ be as in \eqref{7}. Then by the second bullet in \eqref{6}, we may choose a shortest path $W_x$ in $D_x$ from $v_{P_x}$ to some vertex $w_x\in N(x)\cap D_x$ (so $W_x$ is induced and $W_x\setminus w_x$ is disjoint from $N(x)\cap D_x)$.
Let $B_2=Q_2\cup(\bigcup_{x\in D}(W_x\cup P_x))$. If $|D|\geq s$, then by the first and the third bullets of \eqref{6}, $B_2$ contains a subdivision of $T$, and hence so does $G$, a contradiction. So $|D|\leq s-1$. It follows that $Q_2\setminus D$ has at most $s$ connected components, and by \eqref{5}, $Q_2\setminus D$ contains at least $6s^2$ irredundant vertices. Therefore, there exists a connected component $Q_3$ of $Q_2\setminus D$ containing at least $6s$ irredundant vertices. This proves \eqref{8}.\vsp 

From now on, let $Q_3$ be as obtained in \eqref{8},  $r=|Q_3|$, and $Q_3=q_1\dd \cdots \dd q_r$.  By \eqref{8}, we have $r\geq 6s\geq 6$. For every $i\in \{2,\ldots, r-1\}$, we denote by $L_i$ and $R_i$ the components of $Q_2\setminus q_i$ containing $q_{i-1}$ and $q_{i+1}$, respectively. Since $q_i$ is not docile, the vertex-set of every connected component of $H\setminus q_i$ contains either $L_i$ or $R_i$. Also, by \eqref{3}, $q_i$ is a cut-vertex of $H$. So $H\setminus q_i$ has exactly two distinct connected componets $\lambda_i$ and $\rho_i$, such that $L_i\subseteq \lambda_i$ and $R_i\subseteq \rho_i$. For every $i\in \{1,\ldots, r-1\}$, let us say $i$ is a \textit{bump} if there exists a connected component of $H\setminus \{q_i,q_{i+1}\}$, denoted by $\mu_{i}$, such that $\mu_i\cap Q_2=\emptyset$. From this definition, we immediately deduce the following.
	
\sta{\label{9} Let $i \in \{1,\ldots, r-1\}$ be a bump. Then 
\begin{itemize}
		\item $\mu_i$ is anticomplete to $Q_2\setminus \{q_i,q_{i+1}\}$;
		\item $(N(q_{i})\cup N(q_{i+1}))\cap \mu_{i}\neq \emptyset$; and
	    \item for every bump $j\in \{1,\ldots, r-1\}\setminus \{i\}$, $\mu_i$ is anticomplete to $\mu_j$ in $G$.
	\end{itemize}}
	
Also, we have:
	
\sta{\label{10}  For every bump $i\in \{1,\ldots, r-1\}$, there is a (possibly not unique) path $P_i\in \mathcal{P}$ with $v_{P_i}\in \mu_i$.}

Otherwise $H-\mu_i$ and $\mathcal{P}$ satisfy \eqref{12}, violating the minimality of $H$. This proves \eqref{10}.

\sta{\label{11} For every $i\in \{2,\ldots, r-1\}$, if $q_i$ is irredundant, then either $i-1$ or $i$ is a bump.}

Since $q_i$ is irredundant, it has a neighbor $w\in H\setminus \{q_{i-1},q_i,q_{i+1}\}$. It follows that either $w\in \lambda_i\setminus \{q_{i-1}\}$ or $w\in \rho_i\setminus \{q_{i+1}\}$. Assume the former holds. Note that $q_i$ separates $w$ from $R_i$ in $H$. Also, if there exists a path $M$ in $H\setminus q_{i-1}$ from $w$ to some vertex in $x\in L_{i-1}$, then $q_i\dd w\dd M\dd x$ is a path in $H\setminus q_{i-1}$ from $q_i\in R_{i-1}$ to $x\in L_{i-1}$, a contradiction. Thus, $q_{i-1}$ separates $w$ from $L_{i-1}$ in $H$. As a result, $\{q_{i-1},q_i\}$ separates $w$ from $Q_2\setminus \{q_{i-1},q_{i}\}$, and so the connected component of $H\setminus \{q_{i-1},q_{i}\}$ containing $w$ is anticomplete to $Q_2\setminus \{q_{i-1},q_{i}\}$. But then $i-1$ is a bump. Similarly, if $w\in \rho_i\setminus \{q_{i+1}\}$, then $i$ is a bump. This proves \eqref{11}.\vsp

Let $i\in \{1,\ldots,r-1\}$ be a bump, and $P_i\in \mathcal{P}$ be as in \eqref{10}. By the second bullet of \eqref{9}, we may a choose a shortest path $Z_i$ in $\mu_i$ from $v_{P_i}$ to some vertex $z_i\in (N(q_{i})\cup N(q_{i+1}))\cap \mu_{i}$ (so $Z_i$ is induced and $Z_i\setminus z_i$ is disjoint from $ N(q_{i})\cup N(q_{i-1}) $). We say $i$ is a bump of \textit{type 1} if $z_i\in N(q_{i})\setminus  N(q_{i+1})$, of \textit{type 2} if $z_i\in N(q_{i+1})\setminus  N(q_{i})$, and of \textit{type 3} if  $z_i\in N(q_{i})\cap N(q_{i+1})$. Note that every bump is of type 1, 2 or 3.

By \eqref{8}, there exists $I\subseteq \{2,\ldots, r-1\}$ with $|I|\geq 6s-2$ such that $q_i$ is irredundant for all $i\in I$. Therefore, by \eqref{11}, there exists a $I'\subseteq \{1,\ldots, r-1\}$ with $|I'|\geq 3s-1$ such that every $i\in I'$ is a bump. Consequently, there exists $I''\subseteq I$ with $|I''|\geq s$ such that all elements of $I''$ are bumps of the same type. Now, let $B_3=Q_2\cup(\bigcup_{i\in I''}(Z_i\cup P_i))$. If either all elements of $I''$ are of type 1 or all elements of $I''$ are of type 2, then $B_3$ contains a subdivision of $T$, which is impossible. Otherwise, all elements of $I''$ are of type 3. But then $B_3$ contains the line graph of a subdivision of $T$, a contradiction. This proves Theorem \ref{cater}.
\end{proof}
Next we prove a lemma.
\begin{lemma}
  \label{lem:wallandcreature}
  Let $\Delta, b,k,t$ be positive integers where $k \geq 3$.
  Let $\mathcal{C}$  be the class of  graphs with maximum degree $\Delta$
that do not contain a  $(k,t)$-creature or the line graph of a subdivision of $W_{b\times b}$. 
There exists $R_{b, t,k,  \Delta}$ such that $\tw(G) \leq R_{b,t,k,\Delta}$ for
  every $G \in \mathcal{C}$.
\end{lemma}
\begin{proof}

   Let $t_i = t(1 + \Delta)^{k-i}$. Let $\mathcal{C}_i$ be the class of  graphs with maximum degree $\Delta$
that do not contain an  $(i,t_i)$-creature and have no induced subgraph isomorphic to the line graph of a subdivision of $W_{b\times b}$. 
We will prove by induction that there exists $R_{b, t,k, i, \Delta}$ such that $\tw(G) \leq R_{b,t,k,i,\Delta}$ for every $G \in \mathcal{C}_i$. Since $S_{t_3, t_3, t_3}$ is a $(3, t_3)$-creature, for $i=3$ the result follows from  Theorem~\ref{thm:claw_free}.  Next we prove a result about the existence of forcers in graphs in $\mathcal{C}_i$.

  \sta{\label{claws_are_forcers_caterpillar} $S_{t_i+1,t_i+1,t_i+1} + H$ is a $S_{t_i, t_i+1, t_i+1}$-forcer for $\mathcal{C}_i$ for every $(i-1, t_{i-1})$-creature $H$.}

  Let $G \in \mathcal{C}_i$ and let $H$ be an $(i-1, t_{i-1})$-creature. Let $Y$ be an $S_{t_i+1, t_i+1, t_i+1} +H$ in $G$, let  $Y' = Y \setminus H$, let $x \in Y'$ be a leaf of $Y'$, and let $X = Y' \setminus \{x\}$. Let $D$ be a connected component of $G \setminus N[X]$. Suppose $x \in N[D]$. Then, by Lemma \ref{creature}, it follows that $D$ has no $(i-1, t_{i-1})$-creature. Since $H$ is anticomplete to $Y'$, we have that $H \not \subseteq N[D]$. Therefore, $X$ breaks $\{x\} +H$, so $S_{t_i+1, t_i+1, t_i+1} +H$ is a $S_{t_i, t_i+1, t_i+1}$-forcer for $G$. This proves \eqref{claws_are_forcers_caterpillar}.

By Theorem~\ref{thm:centralbag}, it is now enough to bound the treewidth of
 $\{(S_{t_i+1, t_i+1, t_i+1} +H) : H \text{ is an}$ \\ $(i-1, t_{i-1})\text{-creature}\}$-free graphs in $\mathcal{C}_i$. Let $F$ be a graph with no $(i, t_{i})$-creature.
If $F$ is $S_{t_i+1, t_i+1, t_i+1}$-free, the result follows from  Theorem
\ref{thm:claw_free}. Thus, let  $Q \subseteq V(F)$ be  an $S_{t_i+1, t_i+1, t_i+1}$ in $F$. Then, $F \setminus N[Q]$ has no $(i-1, t_{i-1})$-creature, so by the inductive hypothesis, we deduce that $\tw(F \setminus N[Q]) \leq R_{b,t,k,i-1,\Delta}$.
But $|Q|=3t_i+4$, and therefore $|N[Q]| \leq (3t_i+4) \Delta$.
Consequently, $\tw(F) \leq R_{b,t,k,i-1,\Delta}+(3t_i+4) \Delta$, and we can set
$R_{b,t,k,i,\Delta}=R_{b,t,k,i-1,\Delta}+(3t_i+4) \Delta$.
\end{proof}

We can now prove Theorem \ref{thm:caterpillar-non_specific}, which we restate.  

\begin{theorem}
  Let $\Delta$ be a positive integer and let $T$ be a subcubic subdivided caterpillar. Let $\mathcal{C}$  be the class of graphs with maximum degree at most $\Delta$  which do not contain a subdivision of $T$ or the line graph of a subdivision of $T$. Then there exists $R_{\Delta, T}$ such that $\tw(G) \leq R_{\Delta, T}$ for
  every $G \in \mathcal{C}$.
\end{theorem}
\begin{proof}
 
Let $G \in \mathcal{C}$.
  By Theorem \ref{cater}, there exist integers $k,t$ such that if $G \in \mathcal{C}$
 then $G$ does not contain a $(k, t)$-creature. 

 Next we observe:
 
\sta{\label{no_line-graph-wall} Let $G \in \mathcal{C}$. Then $G$ does not contain the line graph of a subdivision of $W_{|T| \times |T|}$.} 

Let $H$ be the line graph of a subdivision of $W_{|T| \times |T|}$. Then, $H$ contains the line graph of a subdivision of $T$. It follows that if $G$ contains $H$, then $G$ contains the line graph of a subdivision of $T$, a contradiction. This proves \eqref{no_line-graph-wall}. \vsp

Now the result follows from Lemma \ref{lem:wallandcreature}.
\end{proof}

\end{document}